\documentclass[a4paper,fleqn]{cas-sc}

\usepackage[authoryear,longnamesfirst]{natbib}

\usepackage{epstopdf}
\usepackage[noabbrev]{cleveref}
\usepackage{mathtools}
\usepackage{mathdots} 
\usepackage{graphicx}
\usepackage{array}              
\usepackage{amsfonts} 
\usepackage{amsmath} 
\usepackage{bm} 
\usepackage{multirow}           
\usepackage{float}               
\usepackage{ulem} 
\usepackage{color}
\usepackage{tikz}
\usetikzlibrary{matrix,shapes,arrows,positioning,chains}

\usepackage{algorithm}
\usepackage{algpseudocode}

\usepackage{subfigure}
\usepackage[shortlabels]{enumitem}
\usepackage{multicol}
\usepackage[makeroom]{cancel}
\usepackage{arydshln}

\newcommand{\llV}{\vspace{-2mm}}
\newcommand{\lllV}{\vspace{-5mm}}
\newcommand{\lH}{\hspace{-1mm}}
\newcommand{\xlH}{\hspace{-.5mm}}
\newcommand{\llH}{\hspace{-2mm}}
\newcommand{\lllH}{\hspace{-5mm}}
\newcommand{\mH}{\hspace{1mm}}
\newcommand{\mmH}{\hspace{2mm}}
\newcommand{\mmmH}{\hspace{5mm}}

\newcommand{\mmV}{\vspace{2mm}}

\newcommand{\N}{\phantom{0}}

\newcommand{\bla}{\color{black}}

\newlength{\Oldarrayrulewidth}

\newcolumntype{"}{@{\hskip\tabcolsep\vrule width 1.5pt\hskip\tabcolsep}}

\newcommand{\T}[1]{\text{#1}} 
\newcommand{\mi}{\T{--}} 
\newcommand{\ma}{\T{+}}

\newcommand{\z}[1]{\bm{#1}}

\newcommand{\br}{\z{r}}
\newcommand{\x}{\z{x}}
\newcommand{\X}{\z{X}}
\newcommand{\A}{\z{A}}
\newcommand{\B}{\z{B}}
\newcommand{\C}{\z{C}}
\newcommand{\D}{\z{D}}
\newcommand{\E}{\z{E}}

\newcommand{\pv}{\bm{\pi}}
\newcommand{\pz}{\pi_0}

\newcommand{\ppk}[1]{\{\{S_k^{#1}\}_{k=1}^p\}_{\succeq}}
\newcommand{\pp}[1]{\bm{\{}S^{#1}\bm{\}}}
\newcommand{\pps}[2]{S_{#2}^{#1}}
\newcommand{\mcc}[2]{\multicolumn{#1}{c}{#2}}
\renewcommand{\N}{\phantom{N}}
\newcommand{\wop}{P_{\T{WO}}^n}
\newcommand{\wopk}[1]{P_{\T{WO}}^{#1}}
\newcommand{\ivl}[3]{#1\xlH\le\xlH #2\xlH\le\xlH #3}

\newcommand{\aw}[2]{\pi_{j_{#1}j_{#2}}} 
\newcommand{\arc}{{\cal A}}
\newcommand{\arcN}{{\cal A}_N}
\newcommand{\arcP}{\arc^{+}}
\newcommand{\arcM}{\arc^{-}}

\newcommand{\fI}{\hat{N}}
\newcommand{\fC}{\hat{N}^c}

\newcommand{\pSize}{\|\x^+\|}
\newcommand{\mSize}{\|\x^-\|}
\newcommand{\lhs}{\pSize\xlH-\xlH\mSize}
\newcommand{\msX}{\z{\bar X}_{M\&R}}

\newcommand{\Tfour}{T1}
\newcommand{\Tfive}{T2-0}
\newcommand{\Tsix}{T2-1}
\newcommand{\Tseven}{T2-2}
\newcommand{\Teight}{T2-4}
\newcommand{\Tnine}{T2-3}

\definecolor{lgray}{gray}{0.8}
\definecolor{dgray}{gray}{0.4}

\newtheorem{example}{Example}

\newtheorem{theorem}{Theorem}
\newtheorem{definition}{Definition}
\newtheorem{lemma}{Lemma}
\newtheorem{conj}{Conjecture}
\newproof{pf}{Proof}
\def\tsc#1{\csdef{#1}{\textsc{\lowercase{#1}}\xspace}}
\tsc{WGM}
\tsc{QE}
\tsc{EP}
\tsc{PMS}
\tsc{BEC}
\tsc{DE}

\begin{document}
\let\WriteBookmarks\relax
\def\floatpagepagefraction{1}
\def\textpagefraction{.001}
\shorttitle{Facet defining inequalities of weak order polytope}
\shortauthors{Adolfo R. Escobedo, Romena Yasmin}

\title [mode = title]{Derivations of large classes of facet defining inequalities of the weak order polytope using ranking structures}    
\tnotemark[1]

\tnotetext[1]{This work was funded by NSF award 1850355.}

\author[1]{Adolfo R. Escobedo}[orcid=0000-0002-4843-3564]
\ead{adRes@asu.edu}
\address[1]{School of Computing, Informatics, and Decision Systems Engineering, Arizona State University, Tempe, AZ}

\author[1]{Romena Yasmin}[orcid=0000-0003-2105-2742]
\cormark[1]
\ead{ryasmin@asu.edu}
\cortext[cor1]{Corresponding author}

\begin{abstract}
Ordering polytopes have been instrumental to the study of combinatorial optimization problems arising in a variety of fields including comparative probability, computational social choice, and group decision-making. The weak order polytope is defined as the convex hull of the characteristic vectors of all binary orders on $n$ alternatives that are reflexive, transitive, and total. By and large, facet defining inequalities (FDIs) of this polytope have been obtained through simple enumeration and through connections with other combinatorial polytopes. This paper derives five new large classes of FDIs by utilizing the equivalent representation of a weak order as a ranking of $n$ alternatives that allows ties; this connection simplifies the construction of valid inequalities, and it enables groupings of characteristic vectors into useful structures. We demonstrate that a number of FDIs previously obtained through enumeration are actually special cases of the large classes. This work also introduces novel construction procedures for generating affinely independent members of the identified ranking structures. Additionally, it states two conjectures on how to derive many more large classes of FDIs using the featured techniques.
\end{abstract}




\maketitle

\section{Introduction}\label{Sec:Intro}
Ordering polytopes have been instrumental to the study of problems from a variety of fields including comparative probability \cite{anscombe1963,chevyrev2013number}, computational social choice \cite{barthelemy1981,marley2016choice} and group decision-making \cite{mar11lin,yoo20new,yoo20cor}. The derivation of facet defining inequalities (FDIs) for these polytopes is a focal point owing primarily to their usefulness in tackling challenging combinatorial optimization problems. For example, their implementation via branch-and-cut approaches has been cited as an effective solution methodology (e.g., see \cite{coll2002facets,gro90fac,mar11lin}). Such polyhedral studies have historically centered on the linear ordering polytope (e.g., see \cite{bolotashvili1999new,fishburn1992induced,grotschel1985facets,marti2011linear}). In recent years, different classes of FDIs have been introduced for various other ordering polytopes (e.g., see \cite{doi01fac,doi16pri,muller1996partial,oswald2003constructing}), driven by a fundamental need to consider a broader range of ordinal relationships. For example, the incorporation of non-strict orderings (i.e., possibly containing ties) is regarded as a fundamental form of preference expression in various group decision-making applications.

There are many situations in which it is necessary to determine an ordering of all alternatives of a set $N:=\{1,\dots,n\}$ or, equivalently a \textit{complete ranking} of $n$ alternatives that best achieves a given objective. For instance, the objective of the consensus ranking problem is to find a ranking $\br\in N^n$ that minimizes the cumulative distance to a set of input rankings $\{\z{a}^k\}_{k=1}^{m}$, the latter of which typically represents a collection of evaluations of the alternatives provided by $m$ sources. It is often prudent, if not necessary, to allow the consensus ranking (as well as the input rankings \cite{kendall1945treatment}) to contain ties for a myriad of reasons---e.g., contradictory information in the given evaluations, high cost of obtaining a strict linear ordering when more than a handful of alternatives are involved, etc. However, the universe of complete non-strict rankings is exceedingly large. This underscores the ongoing need to better characterize the weak order polytope, $\wop$, whose vertices are equivalent to the individual members of this set. Specifically, a \textit{weak order} on $N$ is defined as a binary relation that is reflexive, transitive, and total.

For the most part, FDIs of $\wop$ have been either obtained through enumeration approaches with specialized software or derived from known FDIs of other existing polytopes. The complete sets of FDIs for $\wopk4$ and $\wopk5$ listed in \cite{fio04wea} and \cite{reg12beh}
respectively, were generated using the \texttt{Porta} program \cite{chr97por}. Furthermore, in \cite{doi01fac} certain members of the classes of \textit{2-partition}, \textit{2-chorded cycle}, \textit{2-chorded path}, and \textit{2-chorded wheel} FDIs \citep{gro90fac} of the partial order polytope ${P_{\T{PA}}^n}$ (also known as the clique partitioning polytope) were lifted into FDIs of $\wop$. This work introduces a new approach for deriving large classes of FDIs for $\wop$ from structural insights that are rooted in the connection between weak orders and complete non-strict rankings. 

The main contributions of this work are as follows. We derive six large classes of valid inequalities of $\wop$, for $n\ge4$. Each VI class is obtained by considering a relatively small number of ordinal relationship combinations, specifically the orderings of only up to two alternatives with another and with an unspecified number of other alternatives, which are treated as a group. We prove that five of the six new large classes of VIs are FDIs. This is done by characterizing the weak orders that satisfy a VI at equality into a small number of \textit{ranking structures} and then devising construction procedures that generate expedient sequences of affinely independent points from them. It is important to remark that each fixed value of $n$ provides a VI expression that is not applicable in a dimension lower than $n$; in other words, the statement about the applicability of these FDI classes to any $n\ge4$ is not due to the Lifting Lemma \citep{fio04wea}\bla. In fact, we show that a number of FDIs previously obtained through enumeration are special cases of these large classes of FDIs. Lastly, we conjecture how the above contributions may be extended to derive numerous other large classes of VIs, many of which we expect to be FDIs.

The rest of this paper is organized as follows: \Cref{Sec:Notation} introduces basic  notation and definitions used throughout this work. \Cref{Sec:VI} introduces six classes of valid inequalities of the weak order polytope, along with a description of the ranking structures whose members satisfy each inequality at equality. \Cref{Sec:facets} introduces novel characteristic vector construction procedures, which are then utilized to demonstrate that five of the six aforementioned large classes of valid inequalities are facet defining for any $n\ge4$. \Cref{Sec:Conclusion} concludes with additional insights and two conjectures on how the fundamental insights herein presented may induce additional large classes of FDIs.

\section{Notations, Definitions, and Other Background Concepts} \label{Sec:Notation}
Let $N:=[n]$ be a set of alternatives or objects, where $[n]$ is shorthand for the set $\{1,\dots,n\}$, and let $A_N := \{(i,j): i,j\in N, i\neq j\}$ be the set of all possible ordered pairings on $N$.

Let $\cal W$ be the family of all weak orders on $N$. In certain related works, the  elements of a weak order $W\in\cal W$ are expressed as $i\preceq j$ ($i$ is not preferred over $j$). Without loss of generality, we describe the elements of $W$ as $i\succeq j$ (alternative $i$ is preferred over or tied with $j$) and express this preference relation succinctly as the  ordered pair $(i,j)$. We do this primarily to help make the description of ranking structures more intuitive. Next, we introduce three different representations of $W$ to be used throughout this paper.

\begin{definition}\label{def:charVecRep}
The characteristic vector representation of $W\in\cal W$ is typified by a vector $\x^W\in\{0,1\}^{A_N}$ whose entry $(i,j)\in \arcN$ is given by:
\vspace{-2mm}
\begin{eqnarray*}
  x^W_{(i,j)} = \begin{cases}
  1 & \text{if } (i,j)\in W \\
  0 & \text{otherwise}.
\end{cases}
\end{eqnarray*}
\end{definition}

\begin{definition}\label{def:rankRep}
The (unique) \textit{ranking representation} of $W\in\cal W$ is typified by a vector $\br^W\in N^n$ whose entry $i\in N$ is given by:\llV
\[r^W_i = n - \sum_{j\in N\backslash\{i\}}^{n} \mathbb{I}\left(i\succeq j\in W\right)=n-\sum_{(i,j)\in \arcN} x^W_{(i,j)}.\llV\]
\end{definition}

\begin{definition}\label{def:OrderedPartition}
A collection of nonempty subsets $S_1,\dots,S_p\subseteq N$ forms a \textit{preference partition} of size $p$ of $N$, written as $\ppk{}$, if:
\begin{enumerate}[(a)]
  \item $S_k \cap S_{k'}=\emptyset$ for $k\ne k'$;
  \item $S_1\cup S_2\cup \dots\cup S_p=N$;
  \item $i\approx j$ for $\{i,j\}\in S_k$; $i\succ j$ for $i\in S_k, j\in S_{k'}$ where $1\le k < k'\le p$.
\end{enumerate}
\end{definition}

\begin{definition}\label{def:partitionRep}
The \textit{alternative-ordering representation} of $W\in\cal W$ is typified by a preference partition $\ppk{W}$. The $k$th subset of the partition, written as $S_k^W$, contains all alternatives in the $k$th \textit{bucket} or equivalence class of $W$, where $k\le p$.
\end{definition}

Stated otherwise, $\ppk{W}$ is an ordered set of sets. Its $k$th element, $S_{k}^W$, contains alternatives that are all tied amongst themselves. An example application of the presented terminology is as follows.

\begin{example}\label{ex:Notation}
\normalfont
Let $W=\{(1,2),(1,3),(1,4),(2,1),(2,3),(2,4),(4,3)\}$. The corresponding ranking is $\br^W=(1,1,4,3)$; the rank of alternative 3 is $r^W_3=4$. The corresponding alternative-ordering is ${\{\{S_k^W\}_{k=1}^3\}_{\succeq}}=\{\{1,2\},\{4\},\{3\}\}$; and the contents of bucket 2 are $\pps{W}{2}=\{4\}$.
\end{example}

For notational convenience, when working with a collection of $m$ weak orders, each $W^w\in\cal W$ may be represented in abbreviated fashion as $\x^w, \br^w$, or $\pp{w}$ for $w=1,\dots,m$, when there is no ambiguity; the number of buckets would be evident from the given context or written as $|\pp{w}|$. Additionally, individual elements in the first two representations may be written as $x^w_{ij}$ and $r^w_i$, respectively, for $i,j\in N$ and individual subsets in the third may be written as $\pps{w}{k}$ for $k\le |\pp{w}|$.

Next, we state a few relevant polyhedral theory concepts.  

\begin{definition}\label{def:VI}
The ordered pair $(\pv,\pz)$ is a valid inequality (VI) for a polyhedron $P$ if $\pv\x\le\pz$ holds $\forall \x\in P$ or, equivalently, $\max\{\pv\x:\x\in P\}\le \pi_0$.
\end{definition}

\begin{definition}\label{def:VI_alt}
The valid inequality $(\pv,\pz)$ defines a facet of a polyhedron $P$ if $\emptyset \neq P\cap \{\x:\pv\x=\pz\}\neq P$ and if there exists $\dim(P)$ affinely independent points in $P\cap \{\x:\pv\x=\pz\}$, i.e., the dimension of the VI is $\dim(P)-1$.
\end{definition}

Since $\wop$ has dimension $n(n\mi1)$  (i.e., the polyhedron is full-dimensional), its facets have dimension $n(n\mi1)\mi1$ \cite{gurgel1992}. The FDIs of $\wop$ for $n\geq 3$ that follow directly from the characteristic vector representation and the completeness and transitivity properties of a weak order are \cite{fio04wea}:
\vspace{-2mm}
\begin{subequations}\label{eqn:GKBP}
\begin{eqnarray}
x_{ij} \leq& 1 &\hspace{10mm} i, j = 1,...,n;i \neq j\label{eqn:GKBP1}\\
x_{ij} + x_{ji}  \geq& 1 &\hspace{10mm} i, j = 1,...,n; i \neq j\label{eqn:GKBP2} \\
x_{ij} - x_{ik} - x_{kj}  \geq& -1 &\hspace{10mm} i, j, k = 1,...,n; i \neq j \neq k \neq i \label{eqn:GKBP3}
\end{eqnarray}
\end{subequations}
where, $x_{ij}$ is used in place of $x^W_{(i,j)}$ for visual simplification. 

Later in this work, we refer to the FDIs of $\wopk4$, as categorized into nine classes $WO_1$-$WO_9$ \citep{reg08the} in \cref{WO4_List}. For visual clarity, only nonzero coefficients of $WO_i$ are displayed and the coordinates $(j_k,j_{k'})$ of each VI coefficient-vector $\pv\in\mathbb{R}^{4\times 3}$ are abbreviated as $j_kj_{k'}$, for $(k,k')\in A^{[4]}$.

The cardinality of $WO_i$ in $\wopk4$, written in the rightmost table column as $|WO^4_i|$, is obtained by counting the number of ways that the labels $j_1,j_2,j_3,j_4\in[4]$ can be permuted within each expression. We remark that $WO_1$-$WO_9$ can be reduced into just seven classes \citep{fio01pol,fio04wea} by leveraging symmetries between $WO_6$ and $WO_7$ and between $WO_8$ and $WO_9$. In the list, $WO_1$-$WO_3$  match Inequalities \eqref{eqn:GKBP1}-\eqref{eqn:GKBP3} and are known as the \textit{axiomatic} inequalities, which comprise the full set of FDIs for $\wopk3$. They are also FDIs for $n\ge4$ due to the Lifting Lemma of \citep{fio04wea} restated below for future reference.

\begin{lemma}[Lifting Lemma \citep{fio04wea}]\label{thm:lifting}
Let $(\pv,\pz)$ be an FDI of $\wop$, and let $\bar{\pv}\in\mathbb{R}^{(n+1)\times n}$ be defined by:
\begin{eqnarray*}
  \bar\pi_{jj'} = \begin{cases}
  \pi_{jj'} & \text{ if } j,j'\in[n]: j\ne j',\\
  0 & \text{ if } (j=n\ma1,j'\in[n]) \T{ or } (j\in[n],j'=n\ma1).
\end{cases}
\end{eqnarray*}
Then, $(\bar{\pv},\pz)$ is an FDI of $\wopk{n+1}$.
\end{lemma}

\begin{table}[width=.6\linewidth]\scriptsize
\setlength{\tabcolsep}{2.5pt}
\renewcommand{\arraystretch}{1.5}
\caption{Facet Defining Inequalities $(\pv,\pz)$ of $\wopk4$ and Their Cardinalities}\label{WO4_List}
\begin{tabular}{l|rrrrrrrrrrrr|r|r}$i$ &$\aw12$&$\aw21$&$\aw13$&$\aw31$&$\aw14$&$\aw41$&$\aw23$&$\aw32$&$\aw24$&$\aw42$&$\aw34$&$\aw43$&$\pz$&$|WO^4_i|$\\\hline
1 &1&&&&&&&&&&&&1 &12\\
2 &-1&-1&&&&&&&&&&&-1 &6\\
3 &1&&-1&&&&1&&&&&&-1 &24\\
4 &1    &1  &1  &1  &1  &1  &-1   &-1   &-1   &-1   &-1   &-1   &1   &4\\
5 &-1 &1  &1  &   &1  &   &       &1      &       &1      &-1   &-1   &2   &12\\
6 &1    &   &1  &1  &1  &1  &       &-1   &       &-1   &-1   &-1   &2   &12\\
7 &     &1  &1  &1  &1  &1  &-1   &       &-1   &       &-1   &-1   &2   &12\\
8 &-1 &   &1  &1  &1  &1  &       &1      &       &1      &-1   &-1   &3   &12\\
9 &     &-1&1 &1  &1  &1  &1      &       &1      &       &-1   &-1   &3   &12
\end{tabular}
\end{table}

From the Lifting Lemma, $WO_4$-$WO_9$ are FDIs for any $n\ge4$. Note that the total number of FDIs that can be generated from $WO_4$-$WO_9$ in higher dimensions is greater than in $\wopk4$. Expressly,  the number of these FDIs that can be generated in $\wop$ is given by:
\begin{equation}\label{eqn:WO_Cardinalities}
    |WO^n_i|=\binom{n}{4}|WO^4_i|=\frac{n(n\mi1)(n\mi2)(n\mi3)}{24}\mH|WO^4_i|,
\end{equation}
where $\ivl{4}{i}{9}$ and $n\ge4$. Stated otherwise, $|WO^4_i|$ distinct FDIs from class $WO^4_i$ can be generated for every combination of distinct indices $j_1,j_2,j_3,j_4\in[n]$, each of which can be lifted into any higher dimension. 

\section{Constructing Valid Inequalities}\label{Sec:VI}
This section introduces six new classes of VIs. For fixed dimension $\hat n\ge4$, each class defines a set of VIs specific to $\wopk{\hat n}$, that is, no individual member of the specific VIs for dimension $\hat n$ is applicable in dimension $n < \hat{n}$. While the VIs defined for $\hat n$ are also valid for  $n > \hat n$ owing to the Lifting Lemma, the applicability of the classes to any dimension  $n\ge4$ is independent of lifting. Hence, the cardinality of each VI class is much larger than those of VIs $WO_4$-$WO_9$, which will be elaborated in \cref{Sec:Conclusion}; therein, it will also be explained that these existing FDIs are in fact special cases of the featured VI classes. 

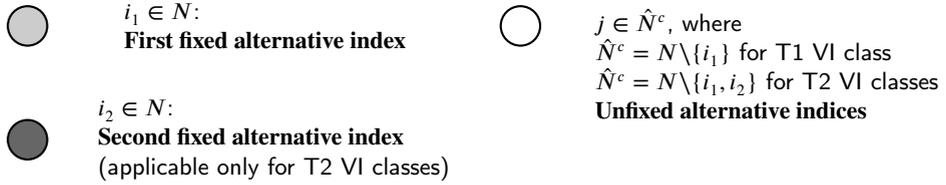
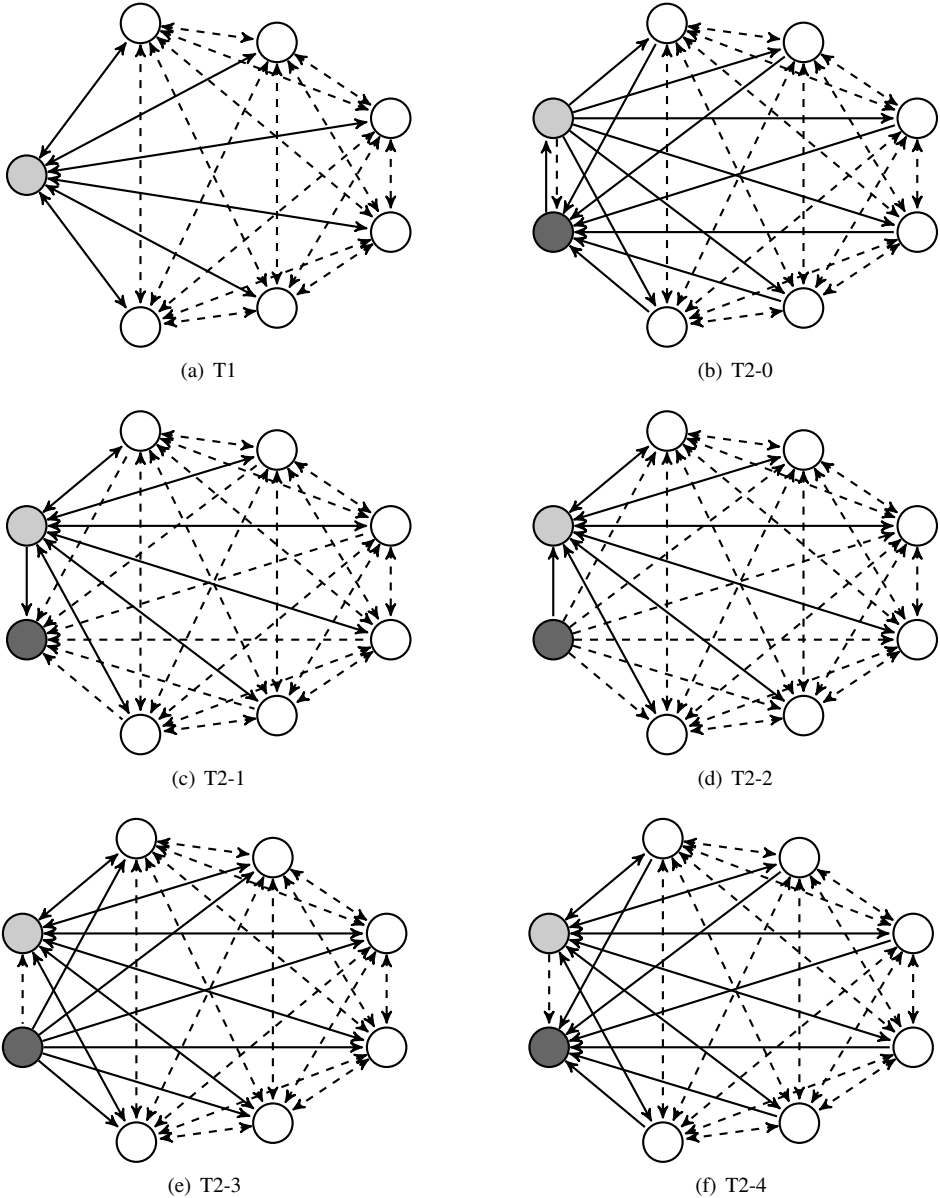
\begin{figure}[!ht]
\centering
\subfigure
{
\begin{tikzpicture}[>=stealth',shorten >=1pt,auto,node distance=2cm, thick,main node/.style={circle,draw,font=\bfseries}]
\node at (0,3.5) [circle,draw,fill=lgray] (i0) {$\phantom{i}$};
\node[align=left] at (3.25,3.5) {$i_1\in N$: \\{\bf First fixed alternative index \phantom{1}}};
\node at (0,2) [circle,draw,fill=dgray] (i1) {$\phantom{i}$};
\node[align=left] at (3.25,2) {$i_2\in N$: \\{\bf Second fixed alternative index}\\ (applicable only for T2 VI classes)};
\node at (6.5,3.5) [circle,draw] (i0) {$\phantom{i}$};
\node[align=left] at (9.75,3) {$j\in\fC$, where\\
$\fC=N\backslash\{i_1\}$ for T1 VI class\\
$\fC=N\backslash\{i_1,i_2\}$ for T2 VI classes\\
\bf Unfixed alternative indices};
\end{tikzpicture}
}
\setcounter{subfigure}{0}

\subfigure[\Tfour] 
{
\begin{tikzpicture}[>=stealth',shorten >=1pt,auto,node distance=2cm, thick,main node/.style={circle,draw,font=\bfseries}]
\node at (0,2) [circle,draw,fill=lgray] (i0) {$\phantom{i}$};
\node at (1.5,4) [circle,draw] (i1) {$\phantom{i}$};
\node at (3.3,3.75) [circle,draw] (i2) {$\phantom{i}$};
\node at (4.8,2.75) [circle,draw] (i3) {$\phantom{i}$};
\node at (4.8,1.25) [circle,draw] (i4) {$\phantom{i}$};
\node at (3.3,0.25) [circle,draw] (i5) {$\phantom{i}$};
\node at (1.5,0) [circle,draw] (i6) {$\phantom{i}$};
\draw [thick,<->] (i0) -- (i1);
\draw [thick,<->] (i0) -- (i2);
\draw [thick,<->] (i0) -- (i3);
\draw [dashed,<->] (i1) -- (i2);
\draw [dashed,<->] (i1) -- (i3);
\draw [thick,<->] (i0) -- (i4);
\draw [thick,<->] (i0) -- (i5);
\draw [thick,<->] (i0) -- (i6);
\draw [dashed,<->] (i4) -- (i5);
\draw [dashed,<->] (i4) -- (i6);
\draw [dashed,<->] (i5) -- (i6);
\draw [dashed,<->] (i1) -- (i4);
\draw [dashed,<->] (i1) -- (i5);
\draw [dashed,<->] (i1) -- (i6);
\draw [dashed,<->] (i2) -- (i4);
\draw [dashed,<->] (i2) -- (i5);
\draw [dashed,<->] (i2) -- (i6);
\draw [dashed,<->] (i3) -- (i4);
\draw [dashed,<->] (i3) -- (i5);
\draw [dashed,<->] (i3) -- (i6);
\draw [dashed,<->] (i2) -- (i3);
\end{tikzpicture}\label{fig:T4}
}\hspace{.5in}
\subfigure[\Tfive]
{
\begin{tikzpicture}[>=stealth',shorten >=1pt,auto,node distance=2cm, thick,main node/.style={circle,draw,font=\bfseries}]
\node at (0,2.75) [circle,draw,fill=lgray] (i0) {$\phantom{i}$};
\node at (1.5,4) [circle,draw] (i1) {$\phantom{i}$};
\node at (3.3,3.75) [circle,draw] (i2) {$\phantom{i}$};
\node at (4.8,2.75) [circle,draw] (i3) {$\phantom{i}$};
\node at (4.8,1.25) [circle,draw] (i4) {$\phantom{i}$};
\node at (3.3,0.25) [circle,draw] (i5) {$\phantom{i}$};
\node at (1.5,0) [circle,draw] (i6) {$\phantom{i}$};
\node at (0,1.25) [circle,draw,fill=dgray] (in) {$\phantom{i}$};
\draw [thick,->] (i0) -- (i1);
\draw [thick,->] (i0) -- (i2);
\draw [thick,->] (i0) -- (i3);
\draw [thick,->] (i0) -- (i4);
\draw [thick,->] (i0) -- (i5);
\draw [thick,->] (i0) -- (i6);
\draw [dashed,->] (i0.280) -- (in.80);
\draw [dashed,<->] (i1) -- (i2);
\draw [dashed,<->] (i1) -- (i3);
\draw [dashed,<->] (i1) -- (i4);
\draw [dashed,<->] (i1) -- (i5);
\draw [dashed,<->] (i1) -- (i6);
\draw [dashed,<->] (i2) -- (i3);
\draw [dashed,<->] (i2) -- (i4);
\draw [dashed,<->] (i2) -- (i5);
\draw [dashed,<->] (i2) -- (i6);
\draw [dashed,<->] (i3) -- (i4);
\draw [dashed,<->] (i3) -- (i5);
\draw [dashed,<->] (i3) -- (i6);
\draw [dashed,<->] (i4) -- (i5);
\draw [dashed,<->] (i4) -- (i6);
\draw [dashed,<->] (i5) -- (i6);
\draw [thick,->] (in.110) -- (i0.250);
\draw [thick,<-] (in) -- (i1);
\draw [thick,<-] (in) -- (i2);
\draw [thick,<-] (in) -- (i3);
\draw [thick,<-] (in) -- (i4);
\draw [thick,<-] (in) -- (i5);
\draw [thick,<-] (in) -- (i6);
\end{tikzpicture}\label{fig:T5}
}

\subfigure[\Tsix]
{
\begin{tikzpicture}[>=stealth',shorten >=1pt,auto,node distance=2cm, thick,main node/.style={circle,draw,font=\bfseries}]
\node at (0,2.75) [circle,draw,fill=lgray] (i0) {$\phantom{i}$};
\node at (1.5,4) [circle,draw] (i1) {$\phantom{i}$};
\node at (3.3,3.75) [circle,draw] (i2) {$\phantom{i}$};
\node at (4.8,2.75) [circle,draw] (i3) {$\phantom{i}$};
\node at (4.8,1.25) [circle,draw] (i4) {$\phantom{i}$};
\node at (3.3,0.25) [circle,draw] (i5) {$\phantom{i}$};
\node at (1.5,0) [circle,draw] (i6) {$\phantom{i}$};
\node at (0,1.25) [circle,draw,fill=dgray] (in) {$\phantom{i}$};
\draw [thick,<->] (i0) -- (i1);
\draw [thick,<->] (i0) -- (i2);
\draw [thick,<->] (i0) -- (i3);
\draw [thick,<->] (i0) -- (i4);
\draw [thick,<->] (i0) -- (i5);
\draw [thick,<->] (i0) -- (i6);
\draw [thick,->] (i0.270) -- (in.90);
\draw [dashed,<->] (i1) -- (i2);
\draw [dashed,<->] (i1) -- (i3);
\draw [dashed,<->] (i1) -- (i4);
\draw [dashed,<->] (i1) -- (i5);
\draw [dashed,<->] (i1) -- (i6);
\draw [dashed,<->] (i2) -- (i3);
\draw [dashed,<->] (i2) -- (i4);
\draw [dashed,<->] (i2) -- (i5);
\draw [dashed,<->] (i2) -- (i6);
\draw [dashed,<->] (i3) -- (i4);
\draw [dashed,<->] (i3) -- (i5);
\draw [dashed,<->] (i3) -- (i6);
\draw [dashed,<->] (i4) -- (i5);
\draw [dashed,<->] (i4) -- (i6);
\draw [dashed,<->] (i5) -- (i6);
\draw [dashed,<-] (in) -- (i1);
\draw [dashed,<-] (in) -- (i2);
\draw [dashed,<-] (in) -- (i3);
\draw [dashed,<-] (in) -- (i4);
\draw [dashed,<-] (in) -- (i5);
\draw [dashed,<-] (in) -- (i6);
\end{tikzpicture}\label{fig:T6}
}\hspace{.5in}
\subfigure[\Tseven]
{
\begin{tikzpicture}[>=stealth',shorten >=1pt,auto,node distance=2cm, thick,main node/.style={circle,draw,font=\bfseries}]
\node at (0,2.75) [circle,draw,fill=lgray] (i0) {$\phantom{i}$};
\node at (1.5,4) [circle,draw] (i1) {$\phantom{i}$};
\node at (3.3,3.75) [circle,draw] (i2) {$\phantom{i}$};
\node at (4.8,2.75) [circle,draw] (i3) {$\phantom{i}$};
\node at (4.8,1.25) [circle,draw] (i4) {$\phantom{i}$};
\node at (3.3,0.25) [circle,draw] (i5) {$\phantom{i}$};
\node at (1.5,0) [circle,draw] (i6) {$\phantom{i}$};
\node at (0,1.25) [circle,draw,fill=dgray] (in) {$\phantom{i}$};
\draw [thick,<->] (i0) -- (i1);
\draw [thick,<->] (i0) -- (i2);
\draw [thick,<->] (i0) -- (i3);
\draw [thick,<->] (i0) -- (i4);
\draw [thick,<->] (i0) -- (i5);
\draw [thick,<->] (i0) -- (i6);
\draw [thick,<-] (i0.270) -- (in.90);
\draw [dashed,<->] (i1) -- (i2);
\draw [dashed,<->] (i1) -- (i3);
\draw [dashed,<->] (i1) -- (i4);
\draw [dashed,<->] (i1) -- (i5);
\draw [dashed,<->] (i1) -- (i6);
\draw [dashed,<->] (i2) -- (i3);
\draw [dashed,<->] (i2) -- (i4);
\draw [dashed,<->] (i2) -- (i5);
\draw [dashed,<->] (i2) -- (i6);
\draw [dashed,<->] (i3) -- (i4);
\draw [dashed,<->] (i3) -- (i5);
\draw [dashed,<->] (i3) -- (i6);
\draw [dashed,<->] (i4) -- (i5);
\draw [dashed,<->] (i4) -- (i6);
\draw [dashed,<->] (i5) -- (i6);
\draw [dashed,->] (in) -- (i1);
\draw [dashed,->] (in) -- (i2);
\draw [dashed,->] (in) -- (i3);
\draw [dashed,->] (in) -- (i4);
\draw [dashed,->] (in) -- (i5);
\draw [dashed,->] (in) -- (i6);
\end{tikzpicture}\label{fig:T7}
}

\subfigure[\Tnine]
{
\begin{tikzpicture}[>=stealth',shorten >=1pt,auto,node distance=2cm, thick,main node/.style={circle,draw,font=\bfseries}]
\node at (0,2.75) [circle,draw,fill=lgray] (i0) {$\phantom{i}$};
\node at (1.5,4) [circle,draw] (i1) {$\phantom{i}$};
\node at (3.3,3.75) [circle,draw] (i2) {$\phantom{i}$};
\node at (4.8,2.75) [circle,draw] (i3) {$\phantom{i}$};
\node at (4.8,1.25) [circle,draw] (i4) {$\phantom{i}$};
\node at (3.3,0.25) [circle,draw] (i5) {$\phantom{i}$};
\node at (1.5,0) [circle,draw] (i6) {$\phantom{i}$};
\node at (0,1.25) [circle,draw,fill=dgray] (in) {$\phantom{i}$};
\draw [thick,<->] (i0) -- (i1);
\draw [thick,<->] (i0) -- (i2);
\draw [thick,<->] (i0) -- (i3);
\draw [thick,<->] (i0) -- (i4);
\draw [thick,<->] (i0) -- (i5);
\draw [thick,<->] (i0) -- (i6);
\draw [thick, dashed,<-] (i0.270) -- (in.90);
\draw [dashed,<->] (i1) -- (i2);
\draw [dashed,<->] (i1) -- (i3);
\draw [dashed,<->] (i1) -- (i4);
\draw [dashed,<->] (i1) -- (i5);
\draw [dashed,<->] (i1) -- (i6);
\draw [dashed,<->] (i2) -- (i3);
\draw [dashed,<->] (i2) -- (i4);
\draw [dashed,<->] (i2) -- (i5);
\draw [dashed,<->] (i2) -- (i6);
\draw [dashed,<->] (i3) -- (i4);
\draw [dashed,<->] (i3) -- (i5);
\draw [dashed,<->] (i3) -- (i6);
\draw [dashed,<->] (i4) -- (i5);
\draw [dashed,<->] (i4) -- (i6);
\draw [dashed,<->] (i5) -- (i6);
\draw [thick,->] (in) -- (i1);
\draw [thick,->] (in) -- (i2);
\draw [thick,->] (in) -- (i3);
\draw [thick,->] (in) -- (i4);
\draw [thick,->] (in) -- (i5);
\draw [thick,->] (in) -- (i6);
\end{tikzpicture}\label{fig:T9}
}\hspace{.5in}
\subfigure[\Teight]
{
\begin{tikzpicture}[>=stealth',shorten >=1pt,auto,node distance=2cm, thick,main node/.style={circle,draw,font=\bfseries}]
\node at (0,2.75) [circle,draw,fill=lgray] (i0) {$\phantom{i}$};
\node at (1.5,4) [circle,draw] (i1) {$\phantom{i}$};
\node at (3.3,3.75) [circle,draw] (i2) {$\phantom{i}$};
\node at (4.8,2.75) [circle,draw] (i3) {$\phantom{i}$};
\node at (4.8,1.25) [circle,draw] (i4) {$\phantom{i}$};
\node at (3.3,0.25) [circle,draw] (i5) {$\phantom{i}$};
\node at (1.5,0) [circle,draw] (i6) {$\phantom{i}$};
\node at (0,1.25) [circle,draw,fill=dgray] (in) {$\phantom{i}$};
\draw [thick,<->] (i0) -- (i1);
\draw [thick,<->] (i0) -- (i2);
\draw [thick,<->] (i0) -- (i3);
\draw [thick,<->] (i0) -- (i4);
\draw [thick,<->] (i0) -- (i5);
\draw [thick,<->] (i0) -- (i6);
\draw [dashed,->] (i0.270) -- (in.90);
\draw [dashed,<->] (i1) -- (i2);
\draw [dashed,<->] (i1) -- (i3);
\draw [dashed,<->] (i1) -- (i4);
\draw [dashed,<->] (i1) -- (i5);
\draw [dashed,<->] (i1) -- (i6);
\draw [dashed,<->] (i2) -- (i3);
\draw [dashed,<->] (i2) -- (i4);
\draw [dashed,<->] (i2) -- (i5);
\draw [dashed,<->] (i2) -- (i6);
\draw [dashed,<->] (i3) -- (i4);
\draw [dashed,<->] (i3) -- (i5);
\draw [dashed,<->] (i3) -- (i6);
\draw [dashed,<->] (i4) -- (i5);
\draw [dashed,<->] (i4) -- (i6);
\draw [dashed,<->] (i5) -- (i6);
\draw [thick,<-] (in) -- (i1);
\draw [thick,<-] (in) -- (i2);
\draw [thick,<-] (in) -- (i3);
\draw [thick,<-] (in) -- (i4);
\draw [thick,<-] (in) -- (i5);
\draw [thick,<-] (in) -- (i6);
\end{tikzpicture}\label{fig:T8}
}
\caption{Digraphs Representing Coefficients $\pv$ of the Six Featured Classes of Valid Inequalities}\label{fig:T4-9_Graphs}
\end{figure}

Surprisingly, each of the featured large classes of VIs can be derived by considering only a small number of ordinal  relationship combinations. \Cref{fig:T4-9_Graphs} helps illustrate this insight with digraphs  $G=(N,\arc)$ that represent the left-hand side coefficients $\pv$ of each of the six classes of VIs, where ${\arc}\subset N\lH\times\lH N$. Each node $j\in N$ represents an alternative and the format of the arc between nodes $j,j'$ represents the coefficient of $x_{jj'}$; expressly, $\pi_{jj'}=0,1,$ or $\mi1$ if there is no arc, a solid arc, or a dashed arc, respectively, from $j$ to  $j'$. When $\pi_{jj'}=\pi_{j'j}\ne0$, a bidirectional arc of the appropriate kind is drawn between  $j$ and $j'$ to represent arcs $(j,j')$ and $(j',j)$. Furthermore, a light gray node represents a ``fixed'' alternative $i_1\in N$; a dark gray node represents a second fixed alternative $i_2\in N$ (when applicable), with $i_2\ne i_1$, and blank nodes represent the remaining ``unfixed'' alternatives $N\backslash\{i_1\}$ or $N\backslash\{i_1,i_2\}$, as applicable. Using this characterization, we categorize the featured VI classes as being of Type 1 or Type 2 (T1 or T2 for short) indicating the number of fixed alternatives in each expression; multiple classes of the same type are differentiated accordingly. Henceforth, the fixed alternative set is denoted as $\fI$, and the unfixed alternative set as $\fC:=N\backslash\fI$ (the complement of $\fI$ in $N$). Additionally, $i$-indices are reserved for elements in $\fI$ while $j$-indices are used for elements of either $\fC$ or $N$, depending on the context. As an important note, although the shaded nodes in each digraph represent the actual number (one or two) of fixed alternatives, the six blank nodes represent a variable number of unfixed alternatives that grows with $n$---more specifically equal to $n-|\fI|$, for $n\ge 4$.

Within each digraph depicted in \ref{fig:T4-9_Graphs}, all arcs belonging to one of the following three sets have a uniform format (all are solid or all are dashed) and a uniform orientation (all point in the same direction or all are bidirectional): $\arc^{\{i_1,j\}}:=\{(i_1,j),(j,i_1): j\in\fC\}, \arc^{\{i_2,j\}}:= \{(i_2,j),(j,i_2): j\in\fC\}$, and $\arc^{\{j,j'\}}:=\{(j,j'): j\in\fC,j\ne j'\}$. In fact, the characteristics of the arcs belonging to the third set are the same across the six digraphs. In other words, even though each of the digraphs (i.e., VI classes) induces a distinctive mathematical expression for any $n\ge4$, their general forms are obtained from a small number of possible uniform-format and uniform-orientation options for $\arc^{\{i_1,j\}}$ and $\arc^{\{i_2,j\}}$, combined with the different arc choices between $i_1$ and $i_2$. For example, in \cref{fig:T6} all arcs in $\arc^{\{i_1,j\}}$ are solid and bidirectional, all arcs in $\arc^{\{i_2,j\}}$ are dashed and directed from $j\in\fC$ to $i_2$, and arc $(i_1,i_2)$ is solid (with no arc from $i_2$ to $i_1$).

The mathematical expressions and derivations of the T1 and T2 VI classes are given in the ensuing paragraphs. In each VI, all non-zero variable coefficients $\pi_{jj'}$ are equal to either 1 or $\mi1$. Hence, the left-hand side of the VI is reexpressed as the difference between the \textit{selected positive arcs} and the \textit{selected negative arcs} from $G$, that is:
\[\lhs:=\sum_{(j,j')\in\arcP}\mH x_{jj'}-\sum_{(j,j')\in\arcM}\mH x_{jj'},\]
where $\arcP:=\{(j,j')\in\arc:\pi_{jj'}=1\}$, $\arcM:=\{(j,j')\in\arc:\pi_{jj'}=-1\}$, $\x^+\in\{0,1\}^{|\arcP|}$, $\x^-\in\{0,1\}^{|\arcM|}$, and $||\cdot||$ is the $L1$-norm. The derivation of each VI centers on determining the value of $\max\{\lhs\}$ in the digraph representations in Figure \ref{fig:T4-9_Graphs}, where $\x^+$ and $\x^-$ together must induce a weak ordering on $N$. Through constructive arguments, the weak orders that achieve this maximum value are described and conveniently grouped into (non-strict) \textit{ranking structures}. The proofs explain why for all possible complete non-strict rankings that do not fit these structures, the value of $\lhs$ must be strictly smaller. For completeness and clarity, the ranking structures are first described in mathematical notation in tables and then in words within the proof narratives. Due to space limitations, some proofs are located in the appendix. Accompanying examples of the weak order constructions described in the VI proofs are given in Supplementary Materials.
 
\begin{theorem}[\Tfour\,VI]\label{theorem:T4 VI}
Let $\fI=\{i_1\}$ be the fixed index set, where $i_1\in N$, and let $\fC=N\backslash\fI$. The following is a valid inequality for ${\bf P}^n_{WO}$, for any $n\ge4$:
\begin{equation}\label{VI_T4}
\sum_{j\in\fC}\left(x_{i_1j}+x_{ji_1}\right)
-\llH\sum_{j,j'\in\fC: j\ne j'}\llH x_{jj'}\le 2 - \frac{(n\mi2)(n\mi3)}{2}.
\end{equation}
\end{theorem}

\begin{table}[width=.6\linewidth,!ht]
\centering
\renewcommand{\arraystretch}{1.25}
\scriptsize
\caption{List of Ranking Structures that Satisfy \Tfour\ VI at Equality}\label{T4_List}
\begin{tabular}{|l|l|c|c|}
  \hline
  \multicolumn{2}{|c|}{Ranking Structure Description}& $\pSize$ & $\mSize$   \\\hline
  \#1& $r_{i_1}=r_{j_1}=k, \T{ where } k\in[n], {j_1}\in\fC$;  & \multirow{2}{*}{$n$} & \multirow{2}{*}{$\frac{(n\mi1)(n\mi2)}{2}$}  \\
  &$r_j, r_{j'}\in[n]\backslash\{k,k\ma1\}, r_j\ne r_{j'}, \T{ for } j,j'\in\fC\backslash\{j_1\}$&&\\\hline 
  \#2& $r_{i_1}=r_{j_1}=r_{j_2}=k,  \T{ where } k\in[n],  {j_1,j_2}\in \fC:j_1\ne j_2$; & \multirow{2}{*}{$n\ma1$} & \multirow{2}{*}{$\frac{(n\mi1)(n\mi2)}{2}\ma1$}  \\
   &$r_j, r_{j'}\in[n]\backslash\{k,k\ma1,k\ma2\}, r_j\ne r_{j'}, \T{ for } j,j'\in\fC\backslash\{j_1,j_2\}$&&\\\hline 
\end{tabular}
\end{table}

\begin{pf}
The VI is satisfied at equality by the ranking structures listed in Table \ref{T4_List}, which are justified as follows. The arc sets $\arcP$ and $\arcM$ are given by:
\begin{flalign*}
 \arcP  & =\{(i_1,j):j\in\fC\} \cup\{(j,i_1):j\in\fC\}, \\
 \arcM  & =\{(j,j'):j,j'\in\fC, j\ne j'\}.
\end{flalign*}
Since there must be at least one arc between every pair of alternatives to induce a complete ranking, by inspection of \cref{fig:T4} the minimum number of selected negative arcs, $\min\mSize$, equals $(n\mi1)(n\mi2)/2$; it is achieved when $j\in\fC$ are ranked in any order but untied. The maximum number of selected positive arcs such that $\mSize$ does not increase can be achieved when $i_1$ is tied with exactly one alternative $j_1\in\fC$, giving precisely $n$ selected positive arcs---2 from the tie of $i_1$ with $j_1$ and $n\mi2$ from the strict ordering of $i_1$ with each $j\in \fC\backslash\{j_1\}$. Hence, $\lhs$ equals the right-hand side of \eqref{VI_T4} for any member of this ranking structure (\#1). This also holds when another positive arc is added to the previous digraph, that is, when $i_1$ is also tied with an alternative $j_2\in \fC\backslash\{j_1\}$, since then $j_1$ and $j_2$ are also tied by transitivity, which requires selecting an additional negative arc (\#2). The value of $\lhs$ achieved by the two structures is maximum since any other ranking structure requires adding at least two negative arcs for each positive arc that is added, owing to transitivity. $_{\square}$\mmV
\end{pf} 

\begin{theorem}[\Tfive\,VI]\label{theorem:T5 VI}
Let $\fI=\{i_1,i_2\}$ be the fixed index set, where $i_1,i_2\in N$ s.t. $i_1\ne i_2$, and let $\fC=N\backslash\fI$. The following is a valid inequality of ${\bf P}^n_{WO}$, for any $n\ge4$:
\begin{equation}\label{VI_T5}
x_{i_2i_1}+\sum_{j\in\fC}\left(x_{i_1j}+x_{ji_2}\right)-x_{i_1i_2} -\llH\sum_{j,j'\in\fC:j\ne j'}\llH x_{jj'}\le
2 - \frac{(n\mi4)(n\mi5)}{2}.
\end{equation}
\end{theorem}\lllV

\begin{table}[width=.6\linewidth,pos=ht]
\centering
\renewcommand{\arraystretch}{1.25}
\scriptsize
\caption{List of Ranking Structures that Satisfy \Tfive\ VI at Equality}\label{T5_List}
\begin{tabular}{|l|l|c|c|}
  \hline
  \multicolumn{2}{|c|}{Ranking Structure Description}& $\pSize$ & $\mSize$   \\\hline
  \#1& $r_{i_1}=1,r_{i_2}=n$; & \multirow{2}{*}{$2(n\mi2)$} & \multirow{2}{*}{$\frac{(n\mi2)(n\mi3)}{2}\ma1$}  \\
  &$r_j, r_{j'}\in[n]\backslash\{1,n\}, r_j\ne r_{j'}, \T{ for } j,j'\in \fC$&&\\\hline 
  \#2& $r_{i_1}=r_{j_1}=1,r_{i_2}=n$, where $j_1\in\fC$; & \multirow{2}{*}{$2(n\mi2)$} & \multirow{2}{*}{$\frac{(n\mi2)(n\mi3)}{2}\ma1$}  \\
  &$r_j, r_{j'}\in[n]\backslash\{1,2,n\}, r_j\ne r_{j'}, \T{ for } j,j'\in \fC\backslash\{j_1\}$&&\\\hline  
   \#3& $r_{i_1}=1,r_{i_2}=r_{j_2}=n\mi1$, where $j_2\in\fC$; & \multirow{2}{*}{$2(n\mi2)$} & \multirow{2}{*}{$\frac{(n\mi2)(n\mi3)}{2}\ma1$}  \\
  &$r_j, r_{j'}\in[n]\backslash\{1,n\mi1,n\}, r_j\ne r_{j'}, \T{ for } j,j'\in \fC\backslash\{j_2\}$&&\\\hline  
   \#4& $r_{i_1}=r_{j_1}=1,r_{i_2}=r_{j_2}=n\mi1$, where $j_1,j_2\in\fC:j_1\ne j_2$; & \multirow{2}{*}{$2(n\mi2)$} & \multirow{2}{*}{$\frac{(n\mi2)(n\mi3)}{2}\ma1$}  \\
  &$r_j, r_{j'}\in[n]\backslash\{1,2,n\mi1,n\}, r_j\ne r_{j'}, \T{ for } j,j'\in \fC\backslash\{j_1,j_2\}$&&\\\hline  
\end{tabular}
\end{table}

\begin{pf}
The VI is satisfied at equality by the ranking structures listed in \cref{T5_List}. Because the VI is only an FDI for $n=4$ (see \cref{Sec:Conclusion}), the remainder of the proof is relegated to \cref{app:T5_VI}.
\end{pf}

\begin{theorem}[\Tsix\,VI]\label{theorem:T6 VI}
Let $\fI=\{i_1,i_2\}$ be the fixed index set, where $i_1,i_2\in N$ s.t. $i_1\ne i_2$, and let $\fC=N\backslash\fI$. The following is a valid inequality of ${\bf P}^n_{WO}$, for any $n\ge4$:
\begin{equation}\label{VI_T6}
x_{i_1i_2}+\sum_{j\in\fC}\left(x_{i_1j}+x_{ji_1}\right)\mmH -\llH\sum_{j,j'\in\fC:j\ne j'}\llH\left(x_{jj'}+x_{ji_2}\right)\le
2 - \frac{(n\mi3)(n\mi4)}{2}.
\end{equation}
\end{theorem}\lllV

\begin{table}[width=.6\linewidth,pos=ht]
\centering
\setlength{\tabcolsep}{2pt}\renewcommand{\arraystretch}{1.25}
\scriptsize
\caption{List of Ranking Structures that Satisfy \Tsix\ VI at Equality}\label{T6_List}
\begin{tabular}{|l|l|c|c|}
  \hline
  \multicolumn{2}{|c|}{Ranking Structure Description}& $\pSize$ & $\mSize$   \\\hline
  \#1& $r_{i_1}=1, r_{i_2}=2;$  & \multirow{2}{*}{$n\mi1$} & \multirow{2}{*}{$\frac{(n\mi2)(n\mi3)}{2}$}  \\
  &$r_j, r_{j'}\in[n]\backslash\{1,2\}, r_j\ne r_{j'}, \T{ for } j,j'\in \fC$&&\\\hline 
  \#2& $r_{i_1}=r_{i_2}=1$; & \multirow{2}{*}{$n\mi1$} & \multirow{2}{*}{$\frac{(n\mi2)(n\mi3)}{2}$}  \\
  &$r_j, r_{j'}\in[n]\backslash\{1,2\}, r_j\ne r_{j'}, \T{ for } j,j'\in \fC$&&\\\hline 
  \#3& $r_{i_1}=r_{j_1}=1, r_{i_2}=3, \T{ where } j_1\in \fC$;  & \multirow{2}{*}{$n$} &  \multirow{2}{*}{$\frac{(n\mi2)(n\mi3)}{2}\ma1$}  \\
  &$r_j, r_{j'}\in[n]\backslash\{1,2,3\}, r_j\ne r_{j'}, \T{ for } j,j'\in \fC\backslash\{j_1\}$&&\\\hline 
  \#4& $r_{i_1}=r_{i_2}=r_{j_1}=1, \T{ where } j_1\in \fC$;  & \multirow{2}{*}{$n$} & \multirow{2}{*}{$\frac{(n\mi2)(n\mi3)}{2}\ma1$}  \\
  &$r_j, r_{j'}\in[n]\backslash\{1,2,3\}, r_j\ne r_{j'}, \T{ for } j,j'\in \fC\backslash\{j_1\}$&&\\\hline 
  \#5& $r_{i_2}=1; r_{i_1}=r_{j_1}=k, \T{ where } k\in[n]\backslash\{1\} \T{ and } j_1\in \fC$;  & \multirow{2}{*}{$n\mi1$} & \multirow{2}{*}{$\frac{(n\mi2)(n\mi3)}{2}$}  \\
  &$r_j, r_{j'}\in[n]\backslash\{1,k,k\ma1\}, r_j\ne r_{j'}, \T{ for } j,j'\in
  \fC\backslash\{j_1\}$&&\\\hline 
  \#6& $r_{i_2}=1, r_{i_1}=r_{j_1}=r_{j_2}=k, \T{ where } k\in[n]\backslash\{1\}, j_1,j_2\in \fC:j_1\ne j_2$;  & \multirow{2}{*}{$n$} & \multirow{2}{*}{$\frac{(n\mi2)(n\mi3)}{2}\ma1$}  \\
  &$r_j, r_{j'}\in[n]\backslash\{1,k,k\ma1,k\ma2\}, r_j\ne r_{j'}, \T{ for } j,j'\in \fC\backslash\{j_1,j_2\}$&&\\\hline 
\end{tabular}
\end{table}

\begin{pf}
The VI is satisfied at equality by the ranking structures listed in Table \ref{T6_List}, which are justified as follows. The arc sets $\arcP$ and $\arcM$ are given by:
\begin{flalign*}
 \arcP  & =\{(i_1,i_2)\}\cup\{(i_1,j):j\in\fC\} \cup\{(j,i_1):j\in\fC\}, \\
 \arcM  & =\{(j,i_2):j\in\fC\}\cup \{(j,j'):j,j'\in\fC, j\ne j'\}.
\end{flalign*}

By inspection of \cref{fig:T6}, $\min\mSize=(n\mi2)(n\mi3)/2$ is achieved when $j\in\fC$ are ranked strictly in any order but all behind $i_2$. The maximum value of $\pSize$ such that $\mSize$ does not increase is equal to $n\mi1$. It is achieved by completing this ranking structure in each of three ways: placing $i_1$ uniquely in first place and $i_2$ uniquely in second place (\#1), tying the two fixed alternatives for first place (\#2), and placing $i_2$ uniquely in first place and tying $i_1$ with some alternative $j_1\in\fC$ for any position inferior to first (\#5). A positive arc can be added to \#1 and \#2 by placing an alternative $j_1\in\fC$ in first place, which ties it with $i_1$ (\#3) or with $i_1$ and $i_2$ (\#4), respectively, adding the negative arc $(j_1,i_2)$ in both cases by transitivity. Similarly, a positive arc can be added to \#5 by tying $i_1$ and $j_1$ with an alternative $j_2\in\fC\backslash\{j_1\}$ for any position inferior to first (\#6), but this also ties $j_1$ and $j_2$ and consequently adds a negative arc. The common value of $\lhs$ achieved by these six ranking structures is maximum. This is because ranking $i_2$ inferior to second place when $i_1$ is uniquely in first place or inferior to first place when $i_1$ and $i_2$ are in a two-way tie only increases  $\mSize$. Moreover, all other ranking structures require adding more positive arcs to one of ranking structures \#1--\#6 (other than those listed), each of which would require adding at least two negative arcs to preserve transitivity. $_{\square}$\mmV
\end{pf}

\begin{theorem}[\Tseven\,VI]\label{theorem:T7 VI}
Let $\fI=\{i_1,i_2\}$ be the fixed index set, where $i_1,i_2\in N$ s.t. $i_1\ne i_2$, and let $\fC=N\backslash\fI$. The following is a valid inequality of ${\bf P}^n_{WO}$, for any $n\ge4$:
\begin{equation}\label{VI_T7}
x_{i_2i_1}+\sum_{j\in\fC}\left(x_{i_1j}+x_{ji_1}\right) -\llH\sum_{j,j'\in\fC:j\ne j'}\llH\left(x_{jj'}+x_{i_2j}\right)\le
2 - \frac{(n\mi3)(n\mi4)}{2}.
\end{equation}
\end{theorem}

\begin{pf}
 The above inequality is satisfied at equality by six ranking structures that are symmetric images of those listed in \cref{T6_List}. That is, those alternatives set to the best and second-best available ranking positions in the \Tsix\, structures are set to the last and second-to-last available positions, respectively, in the \Tseven\, structures. Alternatives occupying the remaining inferior positions in the \Tsix\, structures occupy the remaining superior positions in the \Tseven\, structures. Hence, a string of arguments paralleling the proof of \cref{theorem:T6 VI} establishes that these ranking structures achieve a value of $\max\{\lhs\}$ equal to the right-hand side of \eqref{VI_T7} and all others achieve a lower value. $_{\square}$\mmV
\end{pf} \lllV

\begin{table}[width=.6\linewidth,pos=ht]
\centering
\setlength{\tabcolsep}{2pt}\renewcommand{\arraystretch}{1.25}
\scriptsize
\caption{List of Ranking Structures that Satisfy \Tnine\ VI at Equality}\label{T9_List}
\begin{tabular}{|l|l|c|c|}
  \hline
  \multicolumn{2}{|c|}{Ranking Structure Description}& $\pSize$ & $\mSize$   \\\hline
  \#1& \Tsix\, \#1  & $2(n\mi2)$ & $\frac{(n\mi2)(n\mi3)}{2}$  \\\hline 
  \#2& $r_{i_1}=1, r_{i_2}=r_{j_1}=2, \T{ where } j_1\in \fC$;  & \multirow{2}{*}{$2(n\mi2)$} & \multirow{2}{*}{$\frac{(n\mi2)(n\mi3)}{2}$}  \\
  &$r_j, r_{j'}\in[n]\backslash\{1,2,3\}, r_j\ne r_{j'}, \T{ for } j,j'\in \fC\backslash\{j_1\}$&&\\\hline 
  \#3& \Tsix\, \#3  & $2(n\mi2)$ &  $\frac{(n\mi2)(n\mi3)}{2}$  \\\hline 
  \#4& $r_{i_1}=r_{j_1}=1; r_{i_2}=r_{j_2}=3, \T{ where } j_1,j_2\in\fC:j_1\ne j_2$;  & \multirow{2}{*}{$2(n\mi2)$} & \multirow{2}{*}{$\frac{(n\mi2)(n\mi3)}{2}$}  \\
  &$r_j, r_{j'}\in[n]\backslash\{1,2,3,4\}, r_j\ne r_{j'}, \T{ for } j,j'\in \fC\backslash\{j_1,j_2\}$&&\\\hline
  \#5& \Tsix\, \#4 & $2(n\mi2)\ma1$ & $\frac{(n\mi2)(n\mi3)}{2}\ma1$  \\\hline 
  \#6& $r_{i_1}=r_{i_2}=r_{j_1}=r_{j_2}=1, \T{ where } j_1,j_2\in \fC:j_1\ne j_2$;  & \multirow{2}{*}{$2(n\mi1)$} & \multirow{2}{*}{$\frac{(n\mi2)(n\mi3)}{2}\ma2$}  \\
  &$r_j, r_{j'}\in[n]\backslash\{1,2,3,4\}, r_j\ne r_{j'}, \T{ for } j,j'\in \fC\backslash\{j_1,j_2\}$&&\\\hline 
  \#7& \Tsix\, \#5  & $2(n\mi2)\ma1$ & $\frac{(n\mi2)(n\mi3)}{2}\ma1$ \\\hline 
  \#8& $r_{i_2}=r_{j_2}=1; r_{i_1}=r_{j_1}=k,$ & \multirow{3}{*}{$2(n\mi2)\ma1$} & \multirow{3}{*}{$\frac{(n\mi2)(n\mi3)}{2}\ma1$}  \\
  &$\T{ where }  k\in[n]\backslash\{1,2\}, j_1,j_2\in \fC: j_1\ne j_2$; &&\\
  &$r_j, r_{j'}\in[n]\backslash\{1,2,k,k\ma1\}, r_j\ne r_{j'}, \T{ for } j,j'\in \fC\backslash\{j_1,j_2\}$&&\\\hline  
  \#9& \Tsix\, \#6    & $2(n\mi1)$ & $\frac{(n\mi2)(n\mi3)}{2}\ma2$  \\\hline 
  \#10& $r_{i_2}=r_{j_3}=1; r_{i_1}=r_{j_1}=r_{j_2}=k,$ & \multirow{3}{*}{$2(n\mi1)$} & \multirow{3}{*}{$\frac{(n\mi2)(n\mi3)}{2}\ma2$}  \\
  &$\T{where } k\in[n]\backslash\{1,2\}, j_1,j_2,j_3\in \fC:j_1\ne j_2\ne j_3$;&&\\
   &$r_j, r_{j'}\in[n]\backslash\{1,2,k,k\ma1,k\ma2\}, r_j\ne r_{j'}, \T{ for } j,j'\in \fC\backslash\{j_1,j_2\}$&&\\\hline 
\end{tabular}
\end{table}

\begin{theorem}[\Tnine\,VI]\label{theorem:T9 VI}
Let $\fI=\{i_1,i_2\}$ be the fixed index set, where $i_1,i_2\in N$ s.t. $i_1\ne i_2$, and let $\fC=N\backslash\fI$. The following is a valid inequality of ${\bf P}^n_{WO}$, for any $n\ge4$:
\begin{equation}\label{VI_T9}
\sum_{j\in\fC}\left(x_{i_1j}+x_{ji_1}+x_{i_2j}\right) -x_{i_2 i_1}-\llH\sum_{j,j'\in\fC:j\ne j'}\llH x_{jj'}\le
3 - \frac{(n\mi4)(n\mi5)}{2}.
\end{equation}
\end{theorem}

\begin{pf} The VI is satisfied at equality by the ranking structures listed in \cref{T9_List}. The remainder of the proof is similar to that of \cref{theorem:T6 VI} and can be found in \cref{app:T9_VI}.
\end{pf}

\begin{theorem}[\Teight\,VI]\label{theorem:T8 VI}
Let $\fI=\{i_1,i_2\}$ be the fixed index set, where $i_1,i_2\in N$ s.t. $i_1\ne i_2$, and let $\fC=N\backslash\fI$. The following is a valid inequality of ${\bf P}^n_{WO}$ for any $n\ge4$:
\begin{equation}\label{VI_T8}
\sum_{j\in\fC}\left(x_{i_1j}+x_{ji_1}+x_{ji_2}\right) -x_{i_1 i_2}-\llH\sum_{j,j'\in\fC:j\ne j'}\llH x_{jj'}\le
3 - \frac{(n\mi4)(n\mi5)}{2}.
\end{equation}
\end{theorem}

\begin{pf}
The result can be obtained using a nearly identical line of reasoning as the proof to \cref{theorem:T7 VI}, with the difference that the \Teight\,VI ranking structures are the symmetric images of the \Tnine\,VI ranking structures (see \cref{T9_List}). $_{\square}$\mmV
\end{pf} 

To conclude this section, it is worth restating that \cref{fig:T4-9_Graphs} depicts only a small number of the possible format and orientation choices for $\arc^{\{i_1,j\}}, \arc^{\{i_2,j\}}$, and the arcs that join the fixed alternatives\bla.  \cref{fig:T4} is only 1 of 7  possible digraphs for the case with $|\fI|=1$ and \cref{fig:T5,fig:T6,fig:T7,fig:T8,fig:T9} are only 5 of 343 possible digraphs for the case with $|\fI|=2$. While additional T1 and T2 VIs could be constructed for the omitted combinations using similar constructive arguments, these would likely be equivalent to or dominated by one of the featured VIs. Indeed, as \cref{Sec:Conclusion} demonstrates, \Tfive\, is an FDI only for $n=4$, while \Tfour\, and \Tsix\, to \,\Teight\, are FDIs for all $n\ge4$. It remains an open question whether some of the omitted digraphs represent FDIs for some $n>4$ even though they are not FDIs for $n=4$. Moreover, it remains an open question whether expanding the techniques from this section to cases with $|\fI|\ge3$ (i.e., T3 VIs, T4 VIs, etc.) can produce FDIs of $\wop$ for some or all $n\ge5$. These inquiries are formalized into conjectures in \cref{Sec:Conclusion}.

\section{Constructing Facet Defining Inequalities}\label{Sec:facets} To obtain the dimensionality of the faces induced by the T1 and T2 VIs, we first devise systematic processes for selecting members of their respective ranking structures such that simple patterns of linearly independent characteristic vectors are formed. The key idea is to select these so that consecutive pairs of vectors differ minimally, thereby simplifying the respective proofs.

\subsection{Building Block Procedures} 
Assume that alternatives $i\in I$ belong to bucket $k$ of  alternative-ordering $\pp{w}$, that is, $I\subseteq\pps{w}{k}$, where $1\le k\le p=|\pp{w}|$. Additionally, define a \textit{step parameter} $q\in\mathbb{Q}$, where $\mi k< q < p\mi k\ma1$, $q=t/2$, and $t\in\mathbb{Z}$.
\begin{definition}{}\label{def:push}
A move of $q$ steps of $I$ in $\pp{w}$ is an operation that yields an alternative-ordering in which all $i\in I$ are removed from their current bucket $k$ and either merged with the alternatives in bucket $k\ma q$, when $q\in\mathbb{Z}$, or separated into a new bucket inserted immediately after (before, resp.) bucket $\lfloor k\ma q \rfloor$, when $q\notin\mathbb{Z}$ is positive (negative, resp.). The operation is abbreviated as the triple $\left<I,q,\pp{w}\right>$.
\end{definition}
\begin{example}\label{ex:Move}
\normalfont
Consider five different move operations and their outputs:
\begin{enumerate}
  \item \makebox[5cm]{$\left<\{2\},1,\{\{1,2\},\{4\},\{3\}\}\right>$}$=\{\{1\},\{2,4\},\{3\}\}$
  \item \makebox[5cm]{$\left<\{2,4\},-1,\{\{1\},\{2,4\},\{3\}\}\right>$}$=\{\{1,2,4\},\{3\}\}$
  \item \makebox[5cm]{$\left<\{3\},-2,\{\{1,2\},\{4\},\{3\}\}\right>$}$=\{\{1,2,3\},\{4\}\}$
  \item \makebox[5cm]{$\langle\{1,3\},\frac{3}{2},\{\{1,2,3\},\{4\}\}\rangle$}$=\{\{2\},\{4\},\{1,3\}\}$
  \item \makebox[5cm]{$\langle\{3\},\frac{-5}{2},\{\{1,2\},\{4\},\{3\}\}\rangle$}$=\{\{3\},\{1,2\},\{4\}\}.$
\end{enumerate}
\end{example}
As the example shows, move operations can be used to change not only the contents of a bucket for a given alternative-ordering, but also the ordering and total number of buckets. For instance, the second move operation merges the entire contents of buckets 1 and 2. Additionally, the output alternative-ordering in the third move operation has one fewer bucket than the input alternative-ordering, since the bucket where alternative 3 resides contains only one alternative and $q\in\mathbb{Z}$; the reverse holds for the fourth move operation since the bucket where alternatives 1 and 3 reside contains three alternatives and $q\notin\mathbb{Z}$.

The Merge and Reverse Construction Procedure (M\&R), whose pseudocode is given in \cref{alg:MR}, is at the core of the characteristic vector constructions. It begins with an alternative-ordering $\pp0$ with $p$ buckets (associated with a weak ordering $W^0\in\cal W$), and proceeds to iteratively merge and then reverse adjacent buckets in $\pp0$, generating potential characteristic vectors from each related move operation. To be more precise, the pseudocode displays a shell of the M\&R procedure, which can be customized by incorporating optional steps that allow certain alternatives to move more freely between buckets.
\begin{algorithm}[!ht]
\caption{Merge and Reverse Construction Procedure (M\&R) (Shell)}\label{alg:MR}
\begin{algorithmic}[1]
\Procedure {M\&R}{$\pp0$, $I^0$, $\hat p$}
    \State $p=|\pp0|$ \Comment{number of buckets in input alternative-ordering} 
    \State $I^1=\emptyset$ \Comment{initiate working alternative subset}
    \State $\pp1=\pp0$ \Comment{initiate working alternative-ordering}
    \State $\X=[]$ \Comment{initiate characteristic vector matrix}
    \If{$\hat p\leq p$}
    \For {$j=1,\dots,\hat p$}
        \State $I^1\leftarrow \pps{1}{1}$ \Comment{set $I^1$ to first bucket of $\pp1$}
        \For {$k=1,\dots,p\mi j$}
        \State $\pp1 \leftarrow \left<I^1, 1, \pp1\right>$ \Comment{merge}
        \State $\X$.\textit{append}(\textit{toBinary}$(\pp{1}))$
        \State $\pp1 \leftarrow \left<I^1, \frac{1}{2}, \pp1 \right>$ \Comment{reverse}
        \State $\X$.\textit{append}(\textit{toBinary}$(\pp1))$
        \EndFor
        \State {\bf perform optional outer steps (possibly involving $I^0$)}
    \EndFor
    \EndIf
    \State \textbf{return} $(\X,\pp1)$
\EndProcedure
\end{algorithmic}
\end{algorithm}

\begin{example}\label{ex:HSPC}
\normalfont
Let $\pp0=\{\{1,2\},\{3\},\{4\},\{5\}\}$, $I^0=\{1\}$ and $\hat p = 4$; here, $\pps{0}{1}=\{1,2\}$, $\pps{0}{2}=\{3\},\pps{0}{3}=\{4\},\pps{0}{4}=\{5\}$, and $p=4$. Additionally, define the $j$th optional outer step (pseudocode line 15) for the M\&R as:
\begin{equation}\label{MR_Opt_T4}
\pp1 \leftarrow \left<I^0,j\mi p,\pp1\right>
\end{equation}

Put simply, the optional outer step moves alternative 1 to the first bucket of $\pp1$ after each  execution of the inner for-loop. Performing M\&R($\pp0$, $I^0$, $\hat p$) with this optional outer step produces the following sequence of alternative-orderings:\\
\[\begin{array}{ll||c|c|c}
  j &k  & Merge                     &Reverse                                       &jth\mH Outer\mH Step\\\hline
  1 &1  &\{\{{\bf1,2},3\},\{4\},\{5\}\}  &\{\{3\}\{{\bf1,2}\},\{4\},\{5\}\}&-\\
  1 &2  &\{\{3\}\{{\bf1,2},4\},\{5\}\}   &\{\{3\},\{4\},\{{\bf1,2}\},\{5\}\}&-\\
  1 &3  &\{\{3\},\{4\},\{{\bf1,2},5\}\}  &\{\{3\},\{4\},\{5\},\{{\bf1,2}\}\}& \{\{1,3\},\{4\},\{5\},\{2\}\}  \\
  2 &1  &\{\{{\bf1,3},4\},\{5\},\{2\}\}  &\{\{4\},\{{\bf1,3}\},\{5\},\{2\}\}&-\\
  2 &2  &\{\{4\},\{{\bf1,3},5\},\{2\}\}  &\{\{4\},\{5\},\{{\bf1,3}\},\{2\}\} & \{\{1,4\},\{5\},\{3\},\{2\}\}\\
  3 &1  &\{\{{\bf1,4},5\},\{3\},\{2\}\}  &\{\{5\},\{{\bf1,4}\},\{3\},\{2\}\}&\{\{1,5\},\{4\},\{3\},\{2\}\}\\
\end{array}\]
where the numbers in bold are the members of $I^1$ at iteration $j$. Here, M\&R directly returns 12 characteristic vectors for the weak orders under the ``Merge'' and ``Reverse'' columns, which are then appended to matrix $X$. As a point of emphasis, the optional outer steps serve primarily an auxiliary purpose of setting up $\pp1$ between iterations; characteristic vectors may or may not be stored from their respective weak orders.
\end{example}
\bla
\begin{theorem}\label{thm:MR}
Let $\x^1,\dots,\x^{\hat p(\hat p\mi1)}\in\{0,1\}^{n(n\mi1)}$ denote the $\hat p(\hat p\mi1)$ characteristic vectors generated by the non-optional steps of M\&R($\pp0, I^0, \hat p$), where $p=|\pp0|$ and $\hat p\leq p\leq n$; let $\x^0\in\{0,1\}^{n(n\mi1)}$ denote the corresponding characteristic vector for $\pp0$; and assume that these vectors occupy rows $1,\dots,\hat p(\hat p\mi1)$ and $\hat p(\hat p\mi1)\ma1$, respectively, of $\X\in\{0,1\}^{m\times n(n\mi1)}$, with $m\ge \hat p(\hat p\mi1)\ma1$. Additionally, let $\{j_k\}_{k=1}^{p}$ represent a set of alternative indices, one from each bucket in $\pp0$, which are not permitted to move during the optional outer steps of M\&R; that is, $j_k\in\pps{0}{k}$ with $j_k\notin I^0$, for $k\in[p]$. Independent of the  optional outer steps performed subject to this restriction, M\&R yields at least $\hat p(\hat p\mi1)\ma1$ affinely independent characteristic vectors.
\end{theorem}
\begin{pf}
Assume that $\hat p = p$ in M\&R. The proof focuses on the elements $(j,j')\in \arcN$ such that $j,j'\in\{j_k\}_{k=1}^{p}$. Restricted to these elements, each pair of consecutively generated alternative-orderings differs only in that two alternatives that are in the same bucket in one ordering are in separate adjacent buckets in the other ordering. Stated otherwise, in each successive move operation, $(p\mi2)$ of the alternatives from $\{j_k\}_{k=1}^{p}$ retain their ordinal relationships. 
After subtracting $\x^{i\mi1}$ from $\x^i$, for $i=p(p\mi1),\dots,1$, a submatrix $\msX\in\{0,\mi1,1\}^{p(p\mi1)\ma1\times p(p\mi1)}$ can be extracted from $X$ having the following entry pattern:
\begin{equation}\label{eqn:MR_structure}
\small\msX=\llH \arraycolsep=1pt\def\arraystretch{1}\begin{array}{c|rcrccrc}
\multicolumn{1}{c}{\T{Row}\backslash \T{Col}}&\multicolumn{1}{r}{(j_1,j_2)} & (j_2,j_1) & (j_1,j_3) & (j_3,j_1) &\dots& (j_{p\mi1},j_{p}) & \multicolumn{1}{c}{(j_{ p},j_{p\mi1})} \\\cline{2-8}
1& 0 \mmH\N    &1   &0\mmH\N&0   &\dots                    &0  \mmH\mH\N     &0 \\
2&\mi1 \mmH\N   &0 &0\mmH\N&0   &\dots                    &0  \mmH\mH\N    &0 \\
3&0 \mmH\N   &0 &0\mmH\N&1   &\dots                    &0  \mmH\mH\N    &0 \\
4&0 \mmH\N   &0 &\mi1\mmH\N&0   &\dots                    &0  \mmH\mH\N    &0 \\
\mmH\vdots&\vdots  \mmH\N        &\vdots         &\vdots  \mmH\N       &\vdots             &\ddots        &\vdots\mmH\mH\N & \vdots\\
p(p\mi1)\mi1  &0 \mmH\N   &0 &0\mmH\N&0   &\dots                    &0   \mmH\mH\N   &1 \\
p(p\mi1)\phantom{\mi1}  &0 \mmH\N   &0 &0\mmH\N&0   &\dots                    &\mi1  \mmH\mH\N    &0\\
0& 1 \mmH\N    &0   &1 \mmH\N &0   &\dots                    &1 \mmH\mH\N     &0
\end{array};
\end{equation}
where the 0-row (corresponding to $\x^0$) is placed at the bottom to highlight the convenient structure of this submatrix---expressly, the first $ p(p\mi1)$ rows contain all the possible elementary vectors of size $p(p\mi1)$, times $\pm1$. Next, add the rows with nonzero even indices to row 0 to yield the all-zeros vector, $\z{0}$, of size ${p(p\mi1)}$, and let $\msX'$ be the resulting submatrix. Now, since the rows of the augmented matrix $[\msX'\mH\z{1}]$ are linearly independent, where $\z{1}$ is the all-ones column vector of size $ p(p\mi1)\ma1$, the characteristic vectors $\x^0,\x^1\dots,\x^{ p( p\mi1)}$ are affinely independent. Lastly, when $\hat p < p$,  the resulting matrix $\msX$ that is extracted from the characteristic vectors generated by subroutine M\&R($\pp0, I^0, \hat p$) is contained within the larger matrix obtained with subroutine M\&R($\pp0, I^0, p$) and, thus, the $\hat p(\hat p\mi1)\ma1$ characteristic vectors produced must also be affinely independent. $_{\square}$\mmV
\end{pf} 
 
\subsection{Facet Constructions and Proofs} The presented construction procedures iteratively generate individual characteristic vectors so that they fall into the respective ranking structures of a VI class and differ minimally from preceding vectors (or from other specific reference vectors). These procedures are simplified with the incorporation of M\&R subroutines, which are used to generate a large portion of the needed $n(n\mi1)$ affinely independent vectors; the remaining  vectors are generated to yield other convenient patterns. The respective M\&R subroutines differ in their definition of the optional outer step of the $j$th inner loop; the shorthand statement used to represent a specific M\&R variant in the pseudocode is:
$$\T{M\&R}(S, I, \hat{p})\mH \bm{\vert}\mH j\T{th {\bf optional outer step} defined by (\#)} $$

where (\#) denotes a set of algorithmic expressions. The characteristic vectors are stored in a matrix $\X\in{\mathbb{R}}^{n(n\mi1)\times n(n\mi1)}$. Each row of $\X$ is obtained by converting an alternative-ordering $\{S\}$ into its characteristic vector, which is represented in the pseudocode by the operation $\X.append(toBinary(\{S\}))$. 

Each FDI proof begins with a \textit{difference matrix} $\bar \X\in{\mathbb{R}}^{n(n\mi1)\times n(n\mi1)}$ that reflects the dissimilarities between (mostly) consecutively generated vectors in $\X$; the precise entries of each difference matrix are described in the appendix. Row operations are applied to show that $\bar \X$ is non-singular or, equivalently, that the characteristic vectors are linearly (and affinely) independent. A salient feature of these proofs is that, by leveraging the characteristic vector patterns devised within each construction procedure, proving the non-singularity of $\bar \X$ reduces to showing the non-singularity of a symbolic $4\xlH\times\xlH4$ matrix. For the reader's convenience, numerical examples of the step-by-step matrix operations described within each proof are included in the Supplementary Materials. 

The remainder of this subsection will introduce different construction procedures and demonstrate that five of the six featured VI classes are FDIs for any $n\ge4$. It will also explain why the remaining class is an FDI only for $n=4$. Only one of the FDI proofs is shown within the body of the paper. The others are similar in structure and are located in the appendix due to length considerations. Next, the construction procedure used to show that \Tfour\ VI is an FDI is presented in \cref{alg:T4}; an example of this procedure is also included.   

\begin{example}\label{ex:T4}
\normalfont Perform CP\Tfour ($5,1$) with $j_1=2,j_2=3,j_3=4,j_4=5$ (and $i_1=1$). With these values, line 2 of \cref{alg:T4} initializes $\pp0$ to the starting weak order from \cref{ex:HSPC}, and line 4 yields the 12 characteristic vectors for the weak orders under the ``Merge'' and ``Reverse'' columns therein (i.e., these are the direct outputs of the M\&R subroutine). Line 5 yields the weak order from the last M\&R outer step of \cref{ex:HSPC}---assigned to $\pp1$ in line 4---as the 13th characteristic vector. Through a sequence of move operations that start from $\pp1$, Lines 6-9 yield the next three characteristic vectors and Lines 10-13 three more after that; finally, line 14 yields the initial weak order as the 20th characteristic vector. These last seven weak orders are as follows:\llV
\[\begin{array}{ll}
\{\{4\},\{1,5\},\{3\},\{2\}\}, \\
\{\{4\},\{3\},\{1,5\},\{2\}\}, \\
\{\{4\},\{3\},\{2\},\{1,5\}\}, \\
\{\{4\},\{3\},\{1,2\},\{5\}\}, \\
\{\{4\},\{1,3\},\{2\},\{5\}\}, \\
\{\{1,4\},\{3\},\{2\},\{5\}\},\\
\{\{1,2\},\{3\},\{4\},\{5\}\}.
\end{array}\]
\end{example}

\begin{algorithm}[!t]
\caption{Construction Procedure for Type 1 VI (CP\Tfour)}\label{alg:T4}
\begin{algorithmic}[1]
\Procedure {CP\Tfour}{$n,i_1$}
    \State $\pp0 = \{\{i_1,j_1\},\{j_2\},\{j_3\},\dots,\{j_{n\mi1}\}\}$, where
    $\{j_k\}_{k=1}^{n\mi1}=[n]\backslash\{i_1\}$
    \State $p=|\pp0|=n\mi1$
    \State $(\X,\pp1)\leftarrow $ [ M\&R($\pp0, i_1, p$) $\bm{\vert}$ $j$th {\bf optional outer step} defined by \eqref{MR_Opt_T4} ]
    \State $\X$.\textit{append}(\textit{toBinary}$(\pp1))$
    \For {$j=1,\dots,p\mi1$}
        \State $\pp1\leftarrow$ $\left<\pps{1}{j},\frac{3}{2},\pp1\right>$
        \State $\X$.\textit{append}(\textit{toBinary}$(\pp1))$
    \EndFor
    \For {$j=1,\dots,p\mi1$}
        \State $\pp1\leftarrow$ $\left<i_1,-1,\pp1\right>$
        \State $\X$.\textit{append}(\textit{toBinary}$(\pp1))$
    \EndFor
    \State $\X$.\textit{append}(\textit{toBinary}$(\pp0))$
    \State \textbf{return} $\X$
\EndProcedure
\end{algorithmic}
\end{algorithm}

\begin{theorem}[\Tfour\, FDI]\label{theorem:T4 FDI proof}
\Tfour\,VI is an FDI of ${\bf P}^n_{WO}$, for any $n\ge4$.
\end{theorem}
\begin{pf}
 It is straightforward to verify that each row of $\X$ output by CP\Tfour\, belongs to the ranking structures listed in \cref{T4_List}. 

For ease of exposition, fix $i_1=1$ and $j_k=k\ma1$, for $k=1,\dots,n\mi1$ (or assume a corresponding relabeling of the alternatives is performed a priori). To begin, set $\bar \X$ as the matrix obtained after iteratively subtracting row $i\mi1$ from row $i$, for $i=n(n\mi1)\mi1,\dots,2$, and also subtracting row $i=n(n\mi1)$ from row 1 of $\X$; see \cref{app:T4_Diff} for a full characterization of $\bar \X$. To proceed with row operations, define $\A^0\in{\mathbb{R}}^{n(n\mi1)\times n(n\mi1)}$ and set this matrix with the elements of $\bar{X}$, such that, all comparisons between $j,j'\in\fC=N\backslash\{1\}$ appear in the first $(n\mi1)(n\mi2)$ columns, and the elements involving the comparison of alternative $1$ with $j\in\fC$ show up in the last $2(n\mi1)$ columns. The columns of $\bar \X$ and $\A^0$ are further organized as in \eqref{eqn:MR_structure}; that is, odd columns follow a lexicographical ordering of the respective arcs and even columns have the reverse indexing of the preceding odd columns. Next, add the first $(n\mi1)(n\mi2)/2$ even-index rows to row $n(n\mi1)$ (i.e., the initial weak order) and partition the resulting matrix, $\A^1$, as follows:
\begin{equation}\label{T4_decomposition}
  \A^1=\left[\begin{array}{cc}
\B^1 & \D^1 \\
\C^1 & \E^1
\end{array}\right]
\end{equation}
where, $\B^1\lH \in \mathbb{Z}^{(n\mi1)(n\mi2)\times (n\mi1)(n\mi2)}\lH, \C^1\lH \in\mathbb{Z}^{(2n\mi2)\times (n\mi1)(n\mi2)}\lH, \D^1\lH \in\mathbb{Z}^{(n\mi1)(n\mi2)\times (2n\mi2)}$, and $\E^1\lH \in\mathbb{Z}^{(2n\mi2)\times (2n\mi2)}$. The pertinent entries of the four submatrices are as follows. First, $\B^1$ consists of all but the last row of the $\msX$ submatrix (Equation \eqref{eqn:MR_structure} with $p=n\mi1$), which implies that $|\det(B^1)|=1$. $\C^1$ is mostly a zero matrix, with the exception of row $i$ whose values under columns $(n\mi i\ma1,n)$ and $(n,n\mi i\ma1)$ are $1$ and $\mi1$, respectively, for $i=2,\dots,n\mi1$. Although  $\D^1$ has a more intricate structure than $\B^1$ and $\C^1$, it is only necessary to know the contents of a subset of rows aligned with those rows in $\B^1$ that will be used to turn $\C^1$ into a zero matrix. Expressly, each nonzero row $i$ of $\C^1$ is eliminated to yield an all-zero matrix $\C^2$ by adding to it the two consecutive elementary vectors from $\B^1$ with the opposite signs under columns $(n\mi i\ma1,n)$ and $(n,n\mi i\ma1)$; call them $k(i)$ and $k(i)\ma1$. Rows $k(i)$ and $k(i)\ma1$ of $\D^1$ have a 1 under column $(n,1)$ and a -1 under column $(1,n)$ respectively with no other nonzero entries, for $i=3,\dots,n\mi1$. For, $i=2$, row $k(i)$ of $\D^1$ has -1, 1, and 1 under columns $(n\mi2,1)$, $(1,n\mi1)$, and $(1,n)$, respectively, and no other nonzero entries; row $k(i)\ma1$ of $\D^1$ has -1 under column $(1,n)$ and no other nonzero entries. Finally, $\E^1$ is comprised primarily of a `wraparound staircase'' structure of nonzeros, illustrated as follows:
\[\E^1=\scriptsize\arraycolsep=0.5pt\def\arraystretch{1}
\left[\begin{array}{lcccccccccccccc}
        \mcc{1}{i}& (1,\xlH2) & (2,\xlH1) & (1,\xlH3) & (3,\xlH1) & (1,\xlH4) & (4,\xlH1) & (1,\xlH5) & \dots & (1,\xlH n\mi2) & (n\mi2,\xlH1)&(1,\xlH n\mi1) & (n\mi1,\xlH1) & (1,\xlH n) & (n,1) \\
        1&  &  &  &  &&&&&&  &  & \mi1 & 1 &  \\
        2&  &  &  &  &&&&&& & \mi1 & 1 &  &  \\
        3&  &  &  &  &&&&  & \mi1 & 1 & &&  & \\
        \vdots &  &  &  &  &  & & \iddots &\iddots&&&&&   & \\
        n\mi3&  & &&&  \mi1 & 1 &  &  &  &  &  &  &  &  \\
        n\mi2&  &  & \mi1 & 1 &  &  &  &  &  & && &  &  \\
        n\mi1& \mi1 & 1 &  &  &  &  &  &  & &  &  &  & &  \\
        n  & 1 &  &  &  &  &  &  &  &&&  &  &  & \mi1 \\
        n\ma1&  & \mi1 & 1 &  &  &  &  &  &  &  &  & \\
        n\ma2&  &  &  & \mi1 & 1 & & & &  &  & \\
        n\ma3&  &  &  &&& \mi1 & 1 & & &  &  &  & \\
        \vdots&  & & &  &  &  &&\ddots& \ddots & &  &  &  & \\
        2n\mi3&  &  &  &  &  &&&  &  & \mi1 & 1 &  &  & \\
        2n\mi2& 1 & 1 &  &  & \mi1 & &\mi2&& \dots &  & \mi(n\mi4) &  & \mi(n\mi3) & \\
\end{array}\right];\]\\
$\E^1$'s rows align with rows $(n\mi1)(n\mi 2)\ma1,\dots,n(n\mi1)$ of $\bar \X$, which correspond to the characteristic vectors generated following the M\&R subroutine. We remark that in addition to having a 1 under each of its first two columns, row $2n\mi2$ of $\E^1$ has a decreasing sequence of consecutive negative integers under columns $(1,j)$, for $j=4,\dots,n$. Based on the above explanations, after eliminating the nonzero elements of $\C^1$ using the aforementioned rows from $\B^1$ (and $\D^1$), $\E^1$ changes into $\E^2$, given by:
\[\lH \E^2=\scriptsize\arraycolsep=0.5pt\def\arraystretch{1}
\left[\begin{array}{lcccccccccccccc}
              \mcc{1}{i}& (1,\xlH2) & (2,\xlH1) & (1,\xlH3) & (3,\xlH1) & (1,\xlH4) & (4,\xlH1) & (1,\xlH5) &  \dots & (1,\xlH n\mi2) & (n\mi2,\xlH1)&(1,\xlH n\mi1) & (n\mi1,\xlH1) & (1,\xlH n) & (n,1) \\
              1&  &  &  &  &&&&&&  &  & \mi1 & 1 &  \\
              2&  &  &  &  &&&&&& \mi1& & 1 &  &  \\
              3&  &  &  &  &&&&  & \mi1 & 1 & && \mi1 & 1 \\
              \vdots &  &  &  &  &  & & \iddots &\iddots&&&&& \vdots  & \vdots\\
              n\mi3&  & &&&  \mi1 & 1 &  &  &  &  &  &  & \mi1 & 1 \\
              n\mi2&  &  & \mi1 & 1 &  &  &  &  &  & && & \mi1 & 1 \\
              n\mi1& \mi1 & 1 &  &  &  &  &  &  & &  &  &  & \mi1 & 1 \\
              n  & 1 &  &  &  &  &  &  &  &&&  &  &  & \mi1 \\
              n\ma1&  & \mi1 & 1 &  &  &  &  &  &  &  &  & \\
              n\ma2&  &  &  & \mi1 & 1 & & &  &  &  & \\
              n\ma3&  &  &  &&& \mi1 & 1 & & &  &  &  & \\
              \vdots&  & & &  &  &  &&\ddots& \ddots & &  &  &  & \\
              2n\mi3&  &  &  &  &  &&&  &  & \mi1 & 1 &  &  & \\
              2n\mi2& 1 & 1 &  &  & \mi1 & &\mi2&& \dots &  & \mi(n\mi4) &  & \mi(n\mi3) & \\
\end{array}\right];\]\\
Now, since $\C^2$ is a zero matrix, $|\det(\A^1)|=|\det(\B^1)\det(\E^2)|=|\det(\E^2)|$ and, therefore, we can operate exclusively on $\E^2$ from this point. First, eliminate the nonzero entries along row $2n\mi2$, one by one, from column $(1,2)$ to column $(n\mi2,1)$ by adding to row $2n\mi2$ a multiple of some row $i$, where $\ivl{3}{i}{2n\mi3}$. Specifically, beginning with row $i=n$, alternate between a row with index $i\le n$ and a row with index $i>n$ to select each succeeding pivot row. Note that in each such elimination step, there is only one pivot element available from the designated row-index subset to eliminate the next nonzero entry from row $2n\mi2$, whose value may have been modified by the previous elimination steps. For $4\xlH\le j\xlH<\xlH n$, the sequence value under column $(1,j)$ remains intact until the nonzero under column $(j\mi1,1)$ (the column immediately to the left of $(1,j)$) is eliminated; the latter has a value equal to $1\mi\sum_{k=1}^{j-3}k$, that is, the sum of the preceding negative integers in the sequence, plus one. Second, subtract row $2n\mi3$ (the penultimate row) from row 2 and add rows 3 to $2n\mi4$ to row $2n\mi3$. Upon completion of these row operations, the resulting matrix $\E^3$ possesses the following form:
\[\lH \E^3=\scriptsize\arraycolsep=0.5pt\def\arraystretch{1}
\left[\begin{array}{lccccccccccccccc}
              \mcc{1}{i}& (1,\xlH2) & (2,\xlH1) & (1,\xlH3) & (3,\xlH1) & (1,\xlH4) & (4,\xlH1) & (1,\xlH5) &  \dots & (1,\xlH n\mi2) & (n\mi2,\xlH1)&(1,\xlH n\mi1) & (n\mi1,\xlH1) & (1,\xlH n) & (n,1) \\
              1&  &  &  &  &&&&&&  &  & \mi1 & 1 &  \\
              2&  &  &  &  &&&&&& & \mi1& 1 &  &  \\
              3&  &  &  &  &&&&  & \mi1 & 1 & && \mi1 & 1\\
              \vdots &  &  &  &  &  & & \iddots &\iddots&&&&& \vdots  & \vdots\\
              n\mi3&  & &&&  \mi1 & 1 &  &  &  &  &  &  & \mi1 & 1 \\
              n\mi2&  &  & \mi1 & 1 &  &  &  &  &  & && & \mi1 & 1 \\
              n\mi1& \mi1 & 1 &  &  &  &  &  &  & &  &  &  & \mi1 & 1 \\
              n  & 1 &  &  &  &  &  &  &  &&&  &  &  & \mi1 \\
              n\ma1&  & \mi1 & 1 &  &  &  &  &  &  &  &  & \\
              n\ma2&  &  &  & \mi1 & 1 & & &  &  &  & \\
              n\ma3&  &  &  &&& \mi1 & 1 & & &  &  &  & \\
              \vdots&  & & &  &  &  &&\ddots& \ddots & &  &  &  & \\
              2n\mi3&  &  &  &  &  &&&  &  & & 1 &  & \alpha & \mi \alpha\mi1\\
              2n\mi2&  &  &  &  &  & &&& &  & \beta &  & \gamma & \alpha\mi \gamma\ma1\\
\end{array}\right];\]
where,
\[\alpha=\mi n\ma3;\mmmH \beta=1-\sum_{k=1}^{n-4}k=\frac{-(n\mi3)^2+n-1}{2};\]
\[\gamma=-(n\mi3)-1-\sum_{l=1}^{n\mi5}\left(1-\sum_{k=1}^{l}k\right)=\frac{(n\mi4)^3-13n+46}{6}.\]
Two further points to remark on the structure of $\E^3$ are that columns $(1,2)$ to $(n\mi2,1)$ do not have nonzeros along rows $1,2,2n\mi3,2n\mi2$ and that the submatrix comprised of these $2n\mi6$ columns together with rows 3 to $2n\mi4$ forms a basis---indeed, beginning with column $(1,n\mi1)$, each column can be iteratively added to the left-adjacent column to yield unique elementary columns. Therefore, the first $2n\mi6$ columns can be used to eliminate the nonzero entries along rows 3 to $2n\mi4$ of columns $(1,n\mi1),(n\mi1,1),(1,n),(n,1)$, without impacting the other four rows. Thus, the task of proving the non-singularity of $\E^3$ reduces to proving the non-singularity of the following $4\xlH\times\xlH4$ submatrix induced from its last four columns and first/last two rows:\bla
\[\left[\begin{array}{cccc}
    \,\,\,0     & \mi1    & \,\,\,1      & 0 \\
    \mi1        & \,\,\,1 & \,\,\,0      & 0 \\
    \,\,\,1     & \,\,\,0 & \,\,\,\alpha & \mH \mi \alpha\mi1\\
    \,\,\,\beta & \,\,\,0 & \,\,\,\gamma & \mH\alpha\mi \gamma\ma1
  \end{array}\right].\]
The symbolic determinant of this matrix is $\alpha^2\mi \alpha\beta\mi 2\alpha\mi \beta\mi 1$, which equals 0 if
\[\alpha=\mi 1\Leftrightarrow\mi n\ma3=\mi 1\Leftrightarrow n=4 \T{, or when}\]
\[\alpha=\beta\mi 1\Leftrightarrow\mi n\ma3=\frac{n^2\mi 7n+10}{2}\mi 1\Leftrightarrow n=\frac{9}{2}\pm \frac{\sqrt{17}}{2}\notin \mathbb{Z}.\]
Hence, for $n\ge5$, the $n(n\mi 1)$ characteristic vectors produced by CP\Tfour\, are linearly independent, which implies they are also affinely independent. Lastly, it is straightforward to verify that setting $n=4$ in the \Tfour\ VI expression yields $WO_4$, which was shown to be  facet-defining in \citep{fio04wea} using the \texttt{Porta} program \cite{chr97por}. Therefore \Tfour\,VI is an FDI for any $n\ge 4$. $_{\square}$\mmV
\end{pf}

Next, we describe why \Tfive\ VI is an FDI for $n=4$ but not $n\ge5$. Since $\wop$ is full dimensional, the rank of the characteristic vectors generated by any construction procedure must be at least $n(n\mi1)\mi1$. Because the maximum number of affinely independent weak orders on $N$ that do not contain ties is approximately half of this value, this means that a significant fraction of the generated vectors must contain ties. It can be seen in \cref{T5_List}, however, that the placement of ties is limited to rankings with a two-way tie with alternative $i_1$ for first position (\#2 \& \#4) and rankings with a two way-tie with alternative $i_2$ for last position (\#3 \& \#4). When $n=4$, enough affinely independent vectors can be drawn from the four \Tfive\ VI ranking structures to form a facet---in fact, it is straightforward to verify that substituting this value of $n$ in \Tfive\,VI expression yields $WO_5$. However, as $n$ becomes larger, they do not add up to a sufficiently large fraction. Through enumeration of the weak orders that satisfy the given ranking structures, we verified for many $n$; that the rank of the resulting matrix is far lower than $n(n\mi1)$; for example, for $n=6$, the rank is 15 (i.e., half the required number). Therefore, \Tfive\,VI is not an FDI for $n\ge5$.

The construction procedure used to generate characteristic vectors for  \Tsix\,VI is given in \cref{alg:T6}, where the $j$th optional outer step for the respective M\&R subroutine is defined as:\llV
\begin{subequations}\label{MR_Opt_T6}
\begin{align}
    &\pp1  \leftarrow\bigl<i_2,\frac{_3}{^2},\pp1\bigr> \label{MR_Opt_T6_1} \\
   &\pp1  \leftarrow\left<i_1,j\mi p,\pp1\right> \label{MR_Opt_T6_2}\\
   & \hspace{2mm}\X\T{.\textit{append}(\textit{toBinary}}(\pp1)). \label{MR_Opt_T6_3}
\end{align}
\end{subequations}

\begin{theorem}[\Tsix\,FDI]\label{theorem:T6 FDI proof}
\Tsix\,VI is an FDI of ${\bf P}^n_{WO}$, for any $n\ge4$.
\end{theorem}

\begin{pf}
It is straightforward to verify that all points output by CP\Tsix\,  correspond to the characteristic vectors associated with the six ranking structures that satisfy inequality \eqref{VI_T6} at equality. The remainder of the proof can be found in \cref{app:T6_Diff_&_proof}.
\end{pf} 

\begin{theorem}[\Tseven\,FDI]\label{theorem:T7 FDI proof}
\Tseven\,VI is an FDI of ${\bf P}^n_{WO}$, for any $n\ge4$.
\end{theorem}
\begin{pf}
It is evident from the proof of \cref{theorem:T7 VI} that the applicable characteristic vectors associated with the ranking structures that satisfy inequality \eqref{VI_T7} at equality are symmetric images of those output by CP\Tsix. Therefore, an almost identical set of arguments used in the proof of \cref{theorem:T6 FDI proof} establishes that inequality \eqref{VI_T7} is an FDI for $n\ge4$. $_{\square}$\mmV
\end{pf} 

A comparison of \cref{T6_List,T9_List} reveals that the only ranking structure from \cref{T6_List} that does not satisfy valid inequality \Tnine\, is \#2, in which two fixed alternatives are tied for the first position and the unfixed alternatives are strictly ordered to occupy positions $2$ to $n\mi1$. The characteristic vector that represents ranking structure \#2 from \cref{T6_List} is generated by line 7 of the pseudocode(CP\Tsix), after which a number of move operations are performed to generate a vector representing structure \#3 from the same table. Hence, certain modifications to lines 7-11 are needed to make the construction procedure applicable to the ranking structures associated with \Tnine\,VI. These modifications are as follows. First, after Line 6, move the second fixed alternative, $i_2$, one step to the right to generate a ranking representative of structure \#2 in \cref{T9_List}. Second, break the tie between alternatives $i_2$ and $j_1$ created by the previous step by first moving $i_2$ half a step to the right and then shifting $i_1$ one step to the right to tie it with $j_1$; this generates a ranking representative of structure \#3 in \cref{T9_List}. The above changes can be encapsulated by replacing lines 7-11 of \cref{alg:T6} with the following algorithmic expressions:
\begin{alignat*}{3}
    &\text{Line }7'&:\hspace{20mm}&\pp0  \leftarrow\bigl<i_2,1,\pp0\bigr> \\
   &\text{Line }8'&:\hspace{20mm}&\hspace{2mm} \X.\T{\textit{append}(\textit{toBinary}}(\pp0)) \\
   &\text{Line }9'&:\hspace{20mm}&\pp0   \leftarrow\bigl<i_2,\frac{_1}{^2},\pp0\bigr> \\
   &\text{Line }10'&:\hspace{20mm}&\pp0\leftarrow\bigl<i_1,1,\pp0\bigr>  \\
   &\text{Line }11'&:\hspace{20mm}&\hspace{2mm}\X.\T{\textit{append}(\textit{toBinary}}(\pp0))
\end{alignat*}\lllV

\begin{algorithm}[!ht]
\caption{Construction Procedure for \Tsix\ VI (CP\Tsix)}\label{alg:T6}
\begin{algorithmic}[1]
\Procedure {CP\Tsix}{$n,i_1, i_2$}
    \State $\pp0 =\{\{i_1\},\{i_2\},\{j_1\},\dots,\{j_{n\mi4}\},\{j_{n\mi2}\},\{j_{n\mi3}\}\}$, where $\{j_k\}_{k=1}^{n\mi2}=[n]\backslash\{i_1,i_2\}$
    \State $\X=[]$
    \State $\X$.\textit{append}(\textit{toBinary}$(\pp0))$
    \State $\pp0\leftarrow$ $\left<\pps{0}{n},-\frac{3}{2},\pp0\right>$
    \State $\X$.\textit{append}(\textit{toBinary}$(\pp0))$
    \State $\pp0\leftarrow$ $\left<i_1,1,\pp0\right>$
    \State $\X$.\textit{append}(\textit{toBinary}$(\pp0))$
    \State $\pp0\leftarrow$ $\left<\pps{0}{1},1,\pp0\right>$
    \State $\pp0\leftarrow$ $\left<i_2,\frac{1}{2},\pp0\right>$
    \State $\X$.\textit{append}(\textit{toBinary}$(\pp0))$

    \State $p=|\pp0|=n\mi1$

     \State $(\X,\pp1) \leftarrow $ [ M\&R($\pp0, i_1, p\mi2$) $\bm{\vert}$ $j$th {\bf optional outer step} defined by (\ref{MR_Opt_T6})]
     \State $\pp1\leftarrow$ $\left<i_2,\mi1,\pp1\right>$
    \State $\X$.\textit{append}(\textit{toBinary}$(\pp1))$
    \State $\pp1\leftarrow$ $\left<i_2,\frac{1}{2},\pp1\right>$
    \For {$j=1,\dots,p\mi1$}
        \State $\pp1\leftarrow$ $\left<\pps{1}{j},\frac{3}{2},\pp1\right>$
        \State $\X$.\textit{append}(\textit{toBinary}$(\pp1))$
    \EndFor
    \State \textbf{return} $\X$
\EndProcedure
\end{algorithmic}
\end{algorithm}

\begin{theorem}[\Tnine\,FDI]\label{theorem:T9 FDI proof}
\Tnine\,VI is an FDI of ${\bf P}^n_{WO}$, for $n\ge4$.
\end{theorem}
\begin{pf}
The proof is similar to that of \cref{theorem:T6 FDI proof} and can be found in \cref{app:T9_proof}.
\end{pf} 

\begin{theorem}[\Teight\,FDI]\label{theorem:T8 FDI proof}
\Teight\,VI is an FDI of ${\bf P}^n_{WO}$, for any $n\ge4$.
\end{theorem}
\begin{pf}
This theorem is proved by using a string of arguments similar to the proof of \cref{theorem:T7 FDI proof}, with the difference that the applicable characteristic vectors are symmetric images of those generated by the algorithm used to establish that inequality \eqref{VI_T9} is facet defining for $n\ge4$. $_{\square}$\mmV
\end{pf} 

\section{Additional Insights and Conjectures}\label{Sec:Conclusion} 
We begin this section by discussing how the featured VIs yield FDIs $WO_4$-$WO_9$ (see \cref{Sec:Notation}) as a special case. By fixing $n=4$ in each expression, it is straightforward to verify that VI \Tfour\ yields FDI $WO_4$, while VIs \Tfive, \Tsix, and \Tseven\ correspond to FDIs $WO_5$, $WO_6$, and $WO_7$, respectively; finally VI \Tnine\ induces FDI $WO_9$ and VI \Teight\ induces FDI $WO_8$. According to this correspondence, $O(n^4)$ FDIs can be generated from this single variant of each featured VI in a higher dimension $n$ (see \eqref{eqn:WO_Cardinalities}). On that note, it is important to elaborate on how the featured VI classes provide a distinctive expression for each setting of $n\ge4$. To see this, we present the expressions for \Tfour\ VI for $n=4,5$, with $i_1=1$ as the fixed alternative:
\vspace{-2mm}
\begin{subequations}
\begin{equation}\label{eqn:T14}
\begin{aligned}[b]
\hspace{-2mm}x_{12}+x_{21}+x_{13}+x_{31}+x_{14}+x_{41}-x_{23}-x_{32}-x_{24}-x_{42}-x_{34}-x_{43}
  \le 1.
\end{aligned}
\end{equation}
\vspace{-5mm}
\begin{equation}\label{eqn:T15}
\begin{aligned}[b]
x_{12}+x_{21}+x_{13}+x_{31}&+x_{14}+x_{41}-x_{23}-x_{32}-x_{24}-x_{42}-x_{34}-x_{43}\\
& +x_{15}+x_{51}-x_{25}-x_{52}-x_{35}-x_{53}-x_{45}-x_{54} \le -1.
\end{aligned}
\end{equation}
\end{subequations}

Although the first 12 terms of the left-hand side of \eqref{eqn:T15} mirror the left-hand side of \eqref{eqn:T14}, neither of the two inequalities dominates the other. Specifically, for a VI $(\bm{\pi}, \pi_0)$ to be dominated by a VI $(\bm{\pi'}, \pi'_0)$, the relationships $\bm{\pi'} \geq \mu\bm{\pi}$ and $\pi'_0 \leq \mu\pi_0$ must hold, for some scalar $\mu>0$, with at least one of the two inequalities holding in the strict sense. For the two above inequalities, we have that $\pi'_0=-1< 1=\pi_0$---choosing $\mu=1$ since the coefficients of the variables that do not involve alternative index $j=5$ already match. However, the left-hand side terms related to  $j=5$ are both  $\ma1$ and $\mi1$ in \eqref{eqn:T15}, while they are all zero-valued in \eqref{eqn:T14}, implying that $\bm{\pi'} \not\geq \mu\bm{\pi}$. In general, as $n$ increases, the right-hand side will get smaller. Conversely, on the left-hand side, the coefficients of some newly added variables will become smaller than 0, and the coefficient of other newly added variables will become larger than 0 (hence, the incomparability will hold for any $\mu>0$). Similar arguments can be applied to all of the T2 VIs. Therefore, for each of these large VI classes, each fixed value $\hat n\ge4$ yields a distinctive VI that does not dominate or is dominated by another member of the same class in a  lower dimension.


Since a distinctive set of FDIs can be generated from VIs \Tfour, \Tsix, \Tseven, \Tnine, and \Teight\ for any fixed dimension $\hat n\ge 4$, the number of FDIs that can be generated from these new classes of FDIs is very large. Specifically, the total number of FDIs that can be generated from VI class of Type $i$ in dimension $n$ is given by:
\begin{equation}\label{eqn:Rank_VI_Cardinality}
    |WO^n_{Ti}|=\sum_{\hat{n}=4}^n|WO^{\hat n}_{Ti}|=\sum_{\hat{n}=4}^n \binom{n}{\hat{n}} \binom{\hat{n}}{i} i!  =n!\sum_{\hat{n}=4}^n \frac{1}{(n-\hat{n})!(\hat n-i)!}.
\end{equation}
That is, each of the five VI classes generates $O(n!)$ FDIs in dimension $n$. \bla

We conclude the paper by stating two conjectures related to how the insights herein presented may induce many more large classes of FDIs.
\begin{conj}
\label{conj:C1}
Let $\fI=\{i_1,i_2, i_3\}$ be the fixed index set, where $i_1,i_2, i_3\in N$ s.t. $i_1\ne i_2\ne i_3$ and let $\fC=N\backslash\fI$. The following are T3 VIs and FDIs of ${\bf P}^n_{WO}$, for any $n\ge5$:
\vspace{-0.5mm}
\begin{subequations}\label{VI_3}
\begin{equation}\label{VI_C1_1}
\begin{aligned}[b]
x_{i_1i_2} + x_{i_3i_2} + \sum_{j\in\fC} ( x_{i_1j} &+x_{ji_1}  + x_{i_3j} + x_{ji_3} ) - x_{i_1i_3} - x_{i_3i_1}  - \llH\sum_{j,j'\in\fC:j\ne j'}\lllH \left(x_{jj'}  + x_{ji_2}\right) \le 4 - \frac{(n\mi5)(n\mi6)}{2}
\end{aligned}
\end{equation}

\begin{equation}\label{VI_C1_2}
\begin{aligned}[b]
x_{i_2i_1} + x_{i_2i_3} + \sum_{j\in\fC}(x_{i_1j} & +x_{ji_1} + x_{i_3j} + x_{ji_3}) - x_{i_1i_3} - x_{i_3i_1} - \llH\sum_{j,j'\in\fC:j\ne j'}\lllH \left(x_{jj'}  + x_{i_2j}\right) \le 4 - \frac{(n\mi5)(n\mi6)}{2}
\end{aligned}
\end{equation}

\begin{equation}\label{VI_C1_3}
\begin{aligned}[b]
x_{i_1i_3} + x_{i_2i_3} + \sum_{j\in\fC} (x_{i_1j} & + x_{ji_1}  + x_{i_2j} + x_{ji_2}) - x_{i_1i_2} - x_{i_2i_1}  - \llH\sum_{j,j'\in\fC:j\ne j'}\lllH \left(x_{jj'}  + x_{ji_3}\right) \le 4 - \frac{(n\mi5)(n\mi6)}{2}
\end{aligned}
\end{equation}

\begin{equation}\label{VI_C1_4}
\begin{aligned}[b]
x_{i_3i_1} + x_{i_3i_2} + \sum_{j\in\fC} (x_{i_1j} & +x_{ji_1} + x_{i_2j} + x_{ji_2}) - x_{i_1i_2} - x_{i_2i_1} - \llH\sum_{j,j'\in\fC:j\ne j'}\lllH \left(x_{jj'}  + x_{i_3j}\right) \le 4 - \frac{(n\mi5)(n\mi6)}{2}
\end{aligned}
\end{equation}

\begin{equation}\label{VI_C1_5}
\begin{aligned}[b]
\sum_{j\in\fC}\left(x_{i_1j}+x_{ji_1} +x_{ji_2} + x_{i_3j} + x_{ji_3}\right) & - x_{i_1i_2} - x_{i_1i_3} - x_{i_3i_1} - x_{i_3i_2} - \llH\sum_{j,j'\in\fC:j\ne j'}\lllH x_{jj'} \le 5 - \frac{(n\mi6)(n\mi7)}{2}
\end{aligned}
\end{equation}

\begin{equation}\label{VI_C1_6}
\begin{aligned}[b]
\sum_{j\in\fC}\left(x_{i_1j}+x_{ji_1} +x_{i_2j} + x_{i_3j} + x_{ji_3}\right) & - x_{i_2i_1} - x_{i_1i_3} - x_{i_3i_1} - x_{i_2i_3} - \llH\sum_{j,j'\in\fC:j\ne j'}\lllH x_{jj'} \le 5 - \frac{(n\mi6)(n\mi7)}{2}
\end{aligned}
\end{equation}

\begin{equation}\label{VI_C1_7}
\begin{aligned}[b]
\sum_{j\in\fC}(x_{i_1j}+x_{ji_1} + x_{i_2j} + x_{ji_2} & + x_{ji_3}) - x_{i_1i_2} - x_{i_2i_1} - x_{i_1i_3} - x_{i_2i_3} - \llH\sum_{j,j'\in\fC:j\ne j'}\lllH x_{jj'} \le 5 - \frac{(n\mi6)(n\mi7)}{2}
\end{aligned}
\end{equation}

\begin{equation}\label{VI_C1_8}
\begin{aligned}[b]
\sum_{j\in\fC}(x_{i_1j}+x_{ji_1} + x_{i_2j} + x_{ji_2} & + x_{i_3j}) - x_{i_1i_2} - x_{i_2i_1} - x_{i_3i_1} - x_{i_3i_2} - \llH\sum_{j,j'\in\fC:j\ne j'}\lllH x_{jj'} \le 5 - \frac{(n\mi6)(n\mi7)}{2}
\end{aligned}
\end{equation}
\end{subequations}
\end{conj}
\vspace{5mm}
This conjecture was conceived as a two-step process. In the first step, the following binary programming problem was solved for all possible format and orientation choices of arcs that join only the fixed alternatives (i.e., $\{i_1,i_2,i_3\}=\fI$) combined with the possible format and orientation choices of arc sets $\arc^{\{i_1,j\}}, \arc^{\{i_2,j\}}, \arc^{\{i_3,j\}}$ (where $j\in \fC$ is an unfixed alternative index):
\begin{align*}
    \text{max  } \sum_{(i,j)\in\arcP}\mH x_{ij}-\sum_{(i,j)\in\arcM}\mH x_{ij} &\\
    \text{st.  }\hspace{5mm} \text{Constraints }\eqref{eqn:GKBP1}-\eqref{eqn:GKBP3} & \\
                        x_{ij}  \in \{0,1\} & \hspace{10mm} i,j\in N; i\neq j.
\end{align*}
These format and orientation choices determine the contents of the sets $\arcP$ and $\arcM$; we remark that the characteristics of the arcs that join only the unfixed alternatives, namely $\arc^{\{j,j'\}}:=\{(j,j'): j\in\fC,j\ne j'\}$, were set to the constant format and orientation as the six digraphs given in \cref{fig:T4-9_Graphs} (all dashed and bidirectional). In general, the above binary program can be defined for a specific number of fixed alternatives and solved for various values of $n$ so as to deduce a generalizable expression. The objective function provides an algebraic expression for the left-hand side of a VI, while the sequence of numerical objective function values obtained from solving for consecutive values of $n$ are analyzed to derive an algebraic expression for its right-hand side. 

The second step to conceive the above conjecture involved verifying that the VIs generated in the first step are FDIs for at least a few values of $n$. This was done by enumerating all characteristic vectors that satisfy each inequality in \eqref{VI_3} at equality and verifying that these vectors induce a matrix of rank $n(n\mi1)$, for $5\leq n \leq 10$. Hence, the above conjecture is true for at least these values. It remains to show that the result will hold for any $n\ge11$. Proving this using the techniques introduced in this work would require introducing new construction procedures for each VI, followed by additional formal proofs (which would likely be even more protracted). Hence, it is left for future work. 
\begin{conj}
\label{conj:C2}
Let, $\fI=\{i_1,i_2....i_k\}$ be a subset of $N=[n]$, where $k\leq n\mi 2$ and $i_1,i_2....i_k \in N$ s.t. $i_1\ne i_2.... \ne i_k$. By fixing these $|\fI|$ alternatives, new large classes of FDIs that incorporate at least one variable from each pair of alternatives in $N$ can be obtained using the insights presented in this work.
\end{conj}

Verification of \cref{conj:C2} could lead to the development of new techniques for deriving numerous new large classes of FDIs. Each such class would define a set of FDIs specific to $\wopk{\hat n}$, that is, no individual member of the FDIs specific to dimension $\hat n$ would be defined in dimension $q < \hat{n}$. As such, the FDIs would require the involvement of all alternative pairs on $\hat N$. Conversely, when using the Lifting Lemma to generate FDIs of $\wopk{\hat n}$ from FDIs of $P_{WO}^{q}$, only the alternative pairs from a smaller set are utilized---for example, from $\hat N\backslash \{q\ma1,q\ma2,\dots,\hat n\}$.

\appendix
\section{Valid Inequality Proofs}

\subsection{\Tfive\,VI Proof}\label{app:T5_VI}
\begin{theorem}[\Tfive\,VI]
Inequality \eqref{VI_T5} is a VI of ${\bf P}^n_{WO}$, for any $n\ge4$.
\end{theorem}
\begin{pf}
Inequality \eqref{VI_T5} is satisfied at equality by the characteristic vectors corresponding to the ranking structures listed in \cref{T5_List}, which are justified as follows. Here, the positive and negative arc subsets are given by:
\begin{flalign*}
 \arcP  & =\{(i_2,i_1)\}\cup\{(i_1,j):j\in\fC\} \cup\{(j,i_2):j\in\fC\}, \\
 \arcM  & =\{(i_1,i_2)\}\cup \{(j,j'):j,j'\in\fC, j\ne j'\}.
\end{flalign*}
By inspection of \cref{fig:T5}, $\lhs$ equals the right-hand side of Inequality \eqref{VI_T5} for the ranking structure in which $i_1$ is uniquely in first place, $i_2$ is uniquely in last place, and $j\in\fC$ are in any order but untied (\#1). Expressly, this selects all $2(n\mi2)$ positive arcs with a fixed index on one end and an unfixed index on the other, $(n\mi2)(n\mi3)/2$ negative arcs from among only unfixed alternatives (the smallest number possible), and the negative arc $(i_1,i_2)$. The same value for $\lhs$ is also achieved when this ranking structure is slightly modified in each of three ways: an alternative $j_1\in\fC$ is tied with $i_1$ for first place (\#2), an alternative $j_2\in\fC$ is tied with $i_2$ for the last available position, $n\mi1$ in this case (\#3), or both of these ties occur, with $j_1\ne j_2$ (\#4). No positive or negative arcs are added by these three modifications. All other ranking structures achieve a strictly smaller value of $\lhs$, since each additional tie involving either $i_1$ or $i_2$ and an unfixed alternative requires tying at least two unfixed alternatives by transitivity, which increases the magnitude of $\mSize$ while that of $\pSize$ remains unchanged. Furthermore, tying $i_1$ and $i_2$ leads to a net decrease of $(n\mi1)$ positive arcs, since this would also imply that either $i_2$ is top-ranked or $i_1$ is bottom-ranked.
 $_{\square}$\mmV
\end{pf}

\subsection{\Tnine\,VI Proof}\label{app:T9_VI}
\begin{theorem}[\Tnine\,VI]
Inequality \eqref{VI_T9} is a VI of ${\bf P}^n_{WO}$, for any $n\ge4$.
\end{theorem} 
\begin{pf}
Inequality \eqref{VI_T9} is satisfied at equality by the characteristic vectors corresponding to the ranking structures listed in \cref{T9_List}. Here, the positive and negative arc subsets are given by:
\begin{flalign*}
 \arcP  & =\{(i_2,j):j\in\fC\}\cup\{(i_1,j):j\in\fC\} \cup\{(j,i_1):j\in\fC\}, \\
 \arcM  & =\{(i_2,i_1)\}\cup \{(j,j'):j,j'\in\fC, j\ne j'\}.
\end{flalign*}
By inspection of \cref{fig:T9}, the smallest possible value of $\mSize$ is  $(n\mi2)(n\mi3)/2$, and it corresponds to all $j\in\fC$ being strictly ranked and $i_1$ being ranked ahead of $i_2$. The maximum size of $\pSize$ such that the smallest value of $\mSize$ is maintained is equal to $2(n\mi2)$. It is achieved by completing this ranking structure in each of four ways: placing $i_1$ uniquely in first place and either placing $i_2$ uniquely in second place (\#1) or tied with an alternative $j_1\in\fC$ in second place (\#2); or tying $i_1$ with $j_1$ for first place and either placing $i_2$ uniquely in third place (\#3) or tied with an alternative $j_2\in\fC\backslash\{j_1\}$ in third place (\#4). A positive arc can be added to \#2 by tying $i_1$ with $j_1$ for first place, but this also creates a tie with $i_2$ by transitivity and adds the negative arc $(i_2,i_1)$ (\#5). An additional positive arc can be added to \#5 by tying its three top-ranked alternatives with $j_2$, but this also adds a negative arc from the tie between $j_1$ and $j_2$ (\#6). If in \#5, $i_1$ and $j_1$ remain tied for any position but first, which is either occupied by $i_2$ uniquely (\#7) or jointly with $j_2$ (\#8), the respective sizes of $\pSize$ and $\mSize$ are maintained. A positive arc can be added to structure \#7 by tying $i_1$ and $j_1$ with $j_2$ for any position, except first, which is either occupied by $i_1$ uniquely (\#9) or jointly with an alternative $j_3\in \fC\backslash\{j_1,j_2\}$ (\#10). The difference between $\pSize$ and $\mSize$ achieved by these ten ranking structures is maximum. In particular, ranking $i_2$ inferior to the first ranking position when $i_1$ is behind $i_2$ or inferior to the second-best available ranking position when $i_1$ is ahead of $i_2$ only lowers the cardinality of $\pSize$. All other ranking require adding more positive arcs to structures \#1-\#10 than those discussed, each of which would require adding at least two negative arcs to preserve transitivity. $_{\square}$\mmV
\end{pf} 

\section{Facet Defining Inequality Proofs}
\subsection{\Tfour\,FDI Differences Matrix}\label{app:T4_Diff}
To introduce the differences matrix we fixed $i_1=1$ and $j_k=k\ma1$, for $k=1,\dots,n\mi1$; note that the same relabeling was used in the proof. After iteratively subtracting several rows of $\X$ as described in the \Tfour\, FDI proof (see \cref{theorem:T4 FDI proof}), the generating set for $\bar \X$ is given as:
\[\arraycolsep=.75pt\def\arraystretch{0.9}\scriptsize
\begin{array}{|l|l|c||c|c|}\hline
j& \multicolumn{1}{|c|}{I_j:=i\in} &  p,q \T{-ranges}  & (k,\ell)\in \arcN: \bar{X}_{i,(k,\ell)} =1 & (k,\ell)\in \arcN: \bar{X}_{i,(k,\ell)} =\mi1 \\\hline
1&\{ 1\}   &-& -& (3,\{1,2\}) \\
2&\{(2n\mi p)(p\mi1)\mi2p\ma2q\ma2\}&\ivl{1}{p}{n\mi2},\phantom{p}& - & (1,p\ma q\ma1);(p\ma1,p\ma q\ma1)\\
&&\ivl{1}{q}{n\mi p\mi1}&&\\
3&\{(2n\mi p)(p\mi1)\mi2p\ma2q\ma3\}&\ivl{1}{p}{n\mi3},\phantom{p}& (p\ma q\ma2,1);(p\ma q\ma2,p\ma1) & -\\
&&\ivl{1}{q}{n\mi p\mi2}&&\\
4&\{(2n\mi1)p\mi(p\ma1)^2\ma2\}&\ivl{1}{p}{n\mi3}\phantom{,p}& (1,\{\ell\xlH>\xlH p\ma1\}) ; (p\ma3, p\ma2)& (\{p\ma1\} \cup \{k\xlH>\xlH p\ma3\},1)\\
5&\{ (n\mi1)(n\mi2)\ma1\} &-  & (1,n) & (n\mi1,1) \\
6&\{  (n\mi1)(n\mi2)\ma p\ma1\} &\ivl{1}{p}{n\mi2}\phantom{,p}     & (n\mi p,\{1,n\}) & (\{1,n\},n\mi p)\\
7&\{ (n\mi1)(n\mi2)\ma n\}&-   & (1,2) & (n,1)\\
8&\{  (n\mi1)(n\mi2)\ma p\ma n\}&\ivl{1}{p}{n\mi3}\phantom{,p}     & (1,p\ma2) & (p\ma1,1)\\
9&\{n(n\mi1)\} &\ivl{1}{p}{n\mi1}\phantom{,p} & (p,\{\ell\xlH>\xlH p\});(2,1)& -\\\hline 
\end{array}\]

\subsection{\Tsix\,FDI Differences Matrix and Proof}\label{app:T6_Diff_&_proof}

To introduce the differences matrix for the \Tsix\,FDI proof we fixed $i_1=1$, $i_2=n$ and $j_k=k\ma1$, for $k=1,\dots,n\mi2$. After iteratively subtracting several rows of $\X$ as described in the following proof, the generating set for $\bar \X$ is given as:
\[\arraycolsep=.75pt\def\arraystretch{0.9}\scriptsize
\begin{array}{|l|l|c||c|c|}\hline
j& \multicolumn{1}{|c|}{I_j:=i\in} &  p,q \T{-ranges}  & (k,\ell)\in \arcN: \bar{X}_{i,(k,\ell)} =1 & (k,\ell)\in \arcN: \bar{X}_{i,(k,\ell)} =\mi1 \\\hline
1 & \{1\} & - & (2, \{1,n\}) & (n, \{1,2\}) \\
2 & \{ 2n(p\mi1)\mi p^2\ma3 \}  & \ivl{1}{p}{n\mi2}\phantom{,p} & (n,\{1,p\ma1\})  & - \\
3 & \{ 2n(p\mi1)\mi p^2\ma4 \}  & \ivl{1}{p}{n\mi2}\phantom{,p} & - & (\{1,p\ma 1\}, n) \\
4 & \{ (2n\mi p)(p\mi1)\mi p\ma 2q\ma3\}& \ivl{1}{p}{n\mi3},\phantom{p}& (p\ma q\ma1,\{1,p\ma 1\}) & - \\
&&\ivl{1}{q}{n\mi p\mi2}&&\\
5 &\{ (2n\mi p)(p\mi1)\mi p\ma 2q\ma4\}& \ivl{1}{p}{n\mi3},\phantom{p} & - & (\{1,p\ma 1\}, p\ma q\ma 1)\\
&&\ivl{1}{q}{n\mi p\mi2}&&\\
6 & \{ 2np\mi (p\ma 1)^2\ma2 \} & \ivl{1}{p}{n\mi3}\phantom{,p} & (1,\{\ell \xlH>\xlH p\ma1\}); & (\{p\ma1\}\cup\{k \xlH>\xlH p\ma3\} , 1 );\\
& & & (p\ma2, n) & (n,p\ma2)\\
7 & \{ n(n\mi 2)\ma p \}  & \ivl{1}{p}{n\mi3}\phantom{,p} & (n\mi p\mi 1, \{1,n\mi 1\}) & (\{1,n\mi 1\}, n\mi p\mi 1) \\
8 & \{ n(n\mi 1) \mi 2 \}  & - & (n\mi1, n\mi 2) & (n\mi 2, n\mi 1) \\
9 & \{ n(n\mi 1) \mi 1 \}  & - & -  & (n, 1) \\
10 & \{n(n\mi1)\} & \ivl{2}{p}{n\mi2}\phantom{,p} & (1, \{\ell \xlH> 1\xlH\}); (n,\{\ell \xlH<\xlH n\}); & - \\
&&& (p,\{p\xlH< \xlH \ell \xlH<\xlH n\}) & \\
\hline 
\end{array}\]
\newpage 
\begin{theorem}[\Tsix\,FDI]
\Tsix\,VI is an FDI of ${\bf P}^n_{WO}$, for any $n\ge4$.
\end{theorem}

\begin{pf}
In the \Tsix\, valid inequality expression, fix $i_1=1$ and $i_2=n$ for ease of exposition---or assume a corresponding relabeling of the alternatives is performed a priori. It is straightforward to verify that all points output by CP\Tsix\ belong to the six ranking structures given in \cref{T6_List} that satisfy the inequality at equality. 

To obtain the difference matrix, $\bar \X$, first shift rows 1,2 and 3 to the bottom of the matrix such that they become rows $n(n\mi1)\mi2$, $n(n\mi1)\mi1$ and $n(n\mi1)$ respectively; all other rows are shifted upward. Second, iteratively subtract row $i\mi1$ from row $i$, for $i=n(n\mi1)\mi3,\dots,2$. Third, subtract row $n(n\mi1)$ from row 1, row $n(n\mi1)\mi1$ from row $n(n\mi1)\mi2$ and row $n(n\mi1)\mi2$ from row $n(n\mi1)\mi1$. To proceed, set $\A^0\in{\mathbb{R}}^{n(n\mi1)\times n(n\mi1)}$ with a rearranged column ordering of $\bar \X$ such that, all comparisons between the alternatives $j,j'\in\fC=N\backslash\{1,n\}$ appear in the first $(n\mi2)(n\mi3)$ columns, the next $2(n\mi2)$ columns involve the comparisons between $i_2=n$ and $j\in\fC$ and the finally elements involving the comparison of $i_1=1$ with $j\in\fC$ and $i_2=n$ shows up in the last $2(n\mi1)$ columns. The first thing to remark about the structure of $\A^0$ is that for the submatrix involving the first $(n\mi1)(n\mi2)$ columns and the first $(n\mi1)^2\mi1$ rows, nearly all rows have either a 1 or a $\mi1$ as the only nonzero element and the nonzero occurs under a unique column. The only $n\mi2$ rows that do not fit this pattern are rows $2n(i \mi1)\mi i^2\ma2$, for $i=1,\dots,n\mi2$ which have a 1 and a $\mi1$ under columns $(i\ma1,n)$ and $(n,i\ma1)$, respectively. The two consecutive vectors after each of these $n\mi2$ rows have a nonzero element of the opposite sign under the same columns. In particular, row $2n(i \mi1)\mi i^2\ma3$ has a 1 under column $(n,i\ma1)$ and row $2n(i \mi1)\mi i^2\ma4$ has a $\mi1$ under column $(i\ma1, n)$, where $1\leq i\leq n\mi2$. Another thing to note about $\A^0$ is that, the binary values of its final row corresponds to the alternative-ordering in which items $\{1,n\}$ are tied for the first position and the remaining alternatives are in a lexicographical linear ordering occupying positions $2$ to $(n\mi1)$. To eliminate the nonzero elements of the first $(n\mi1)(n\mi2)$ columns of this row add to it rows $2n(i \mi1)\mi i^2\ma2 j$, where, $1\leq i \leq n \mi2$ and $1\leq j\leq n\mi i$. As the next step, eliminate the first two nonzero entries of row $2n(i \mi1)\mi i^2\ma2$ by adding to it rows $2n(i \mi1)\mi i^2\ma3$ and $2n(i \mi1)\mi i^2\ma4$, for $i=1,\dots,n\mi2$. Then shift these $n\mi2$ rows to the bottom of the matrix and denote the resulting matrix as $A^1$. More explicitly, row $(n\mi1)^2\ma i\ma1$ of $\A^1$ receives row $2n(i \mi1)\mi i^2\ma2$ from $\A^0$, for $1\le i\le n\mi2$; all other rows are shifted upwards. Afterwards, the structure of $\A^1$ can be described via a partition with the same number of submatrices and related dimensions as defined by \eqref{T4_decomposition}. 

Similar to the proof of \cref{theorem:T4 FDI proof}, it is only necessary to know a part of the contents of these submatrices to proceed. $\B^1$ is comprised entirely of positive or negative unit vectors and, thus we have, $|\det(\B^1)|=1$. $\C^1$ is mostly a zero matrix, with the exception of row $i$ whose values under columns $(n\mi i\mi1,n\mi1)$ and $(n\mi1,n\mi i\mi1)$ are 1 and \mi1, respectively, for $i=1,..,n\mi3$ and row $n\mi2$ which has a $\mi1$ under column $(n\mi 2,n\mi1)$ and a 1 under column $(n\mi1,n\mi 2)$. To turn $\C^1$ into a zero matrix first add to row $i$, where $1\leq i\leq n\mi3$, the two consecutive elementary vectors from $\B^1$ that have nonzeroes of the opposite sign under columns $(n\mi i\mi1,n\mi1)$ and $(n\mi1,n\mi i\mi1)$. Next, to eliminate the entries of row $n\mi2$ subtract from it rows $(n\mi1)(n\mi2)\mi4$ and $(n\mi1)(n\mi2)\mi3$ of $\B^1$. Similar to $\B^1$, $\D^1$ is comprised entirely of elementary vectors, more specifically, it consists of the following entries: a $1$ under column $(n, 1)$ and $(i\ma j\ma 1,1)$ in row $1$ and $(2n\mi i)(i\mi1)\mi 2(i \mi j)\ma 3$ respectively and a $\mi1$ under $(1,n)$ and $(1, i\ma j\ma 1)$ in row $2$ and $(2n\mi i)(i\mi1)\mi 2(i \mi j)\ma 4$ respectively, where, $1\leq i\leq n\mi3$ and $1\leq j\leq n\mi i\mi 1$. Finally, the structure of $\E^1$ can be described as follows:

\[\lH \E^1=\scriptsize\arraycolsep=.25pt\def\arraystretch{1}\left[\begin{array}{lccccccccccccccc}
              \mcc{1}{i}& (1,\xlH2) & (2,\xlH1) & (1,\xlH3) & (3,\xlH1) & (1,\xlH4) & (4,\xlH1) & \dots & (1,\xlH n\mi3) & (n\mi3,\xlH1)& (1,\xlH n\mi2) & (n\mi2,\xlH1)&(1,\xlH n\mi1) & (n\mi1,\xlH1) & (1,\xlH n) & (n,1) \\
              1&  &  &  &  & & & & & & \mi1 & 1 &  &  &  &  \\
              2&  &  &  &  & & & &\mi1 & 1 &  &  &  &  &  &  \\
              \vdots &  &  &  & &&& \iddots &\iddots&&&& \vdots  & \vdots&&\\
              n\mi4&  &  &\mi1 &1 &  &  &  &  & & & & &  & &  \\
              n\mi3& \mi1 & 1 & & &  &  &  &  & & & & &  & &  \\
              n\mi2&  & & & & & & & & & & & & & & \\
              n\mi1&  &  &  &  &  &  &  &  && &   && & & \mi1\\
              n& 1 &  & 1 & \mi1 & 1 & \mi2 & \dots & 1 & \mi(n\mi5) & 1 & \mi(n\mi4) & 1 & \mi(n\mi4) & &\mi(n\mi3)\\
              n\ma1&  & 1 & & & & & & & & & & & &\mi1 & \\
              n\ma2&  &\mi1 &1 & & 1&\mi1 & \dots&1 & \mi1&1 &\mi1 &1 &\mi1 & & \\
              n\ma3&  & & & \mi1& 1& &\dots &1 & \mi1&1 &\mi1 &1 &\mi1 & & \\
              \vdots&  & & &&& & \ddots & \ddots &&& \vdots  & \vdots & \vdots & \vdots & \vdots \\
              2n\mi3&  &  &  &  &  & && & \mi 1 & 1 &  & 1 & \mi1 &  & \\
              2n\mi2&  &  &  &  &  & && &&& \mi1 & 1 &  &  &
            \end{array}\right];\]
From the above structure we can see that, row $n$ of $\E^1$ in addition to having a 1 under column $(1, i)$, where $1\leq i\leq n\mi1$, also has a non-increasing sequence of negative integers under column $(j, 1)$, where, $3\,\lH \leq j \leq\lH\, n\mi 2$. Upon completion of the elimination steps to convert $\C^1$ into a zero matrix, the entries of $\E^1$ change slightly and only affect the entries in columns $(1,n\mi1)$ and $(n\mi1,1)$. The structure of the new matrix, denoted as $\E^2$, is given by:

\[\lH \E^2=\scriptsize\arraycolsep=.25pt\def\arraystretch{1}\left[\begin{array}{lccccccccccccccc}
              \mcc{1}{i}& (1,\xlH2) & (2,\xlH1) & (1,\xlH3) & (3,\xlH1) & (1,\xlH4) & (4,\xlH1) & \dots & (1,\xlH n\mi3) & (n\mi3,\xlH1)& (1,\xlH n\mi2) & (n\mi2,\xlH1)&(1,\xlH n\mi1) & (n\mi1,\xlH1) & (1,\xlH n) & (n,1) \\
              1&  &  &  &  & & & & & & \mi1 & 1 & \mi1 & 1 &  &  \\
              2&  &  &  &  & & & &\mi1 & 1 &  &  & \mi1 & 1 &  &  \\
              \vdots &  &  &  & &&& \iddots &\iddots&&&& \vdots  & \vdots&&\\
              n\mi4&  &  &\mi1 &1 &  &  &  &  & & & & \mi1 & 1 &  &  \\
              n\mi3& \mi1 & 1 & & &  &  &  &  & & & & \mi1 & 1 &  &  \\
              n\mi2&  & & & & & & & & & & & 1 & \mi1 &  &  \\
              n\mi1&  &  &  &  &  &  &  &  && &   && & & \mi1\\
              n& 1 &  & 1 & \mi1 & 1 & \mi2 & \dots & 1 & \mi(n\mi5) & 1 & \mi(n\mi4) & 1 & \mi(n\mi4) & &\mi(n\mi3)\\
              n\ma1&  & 1 & & & & & & & & & & & &\mi1 & \\
              n\ma2&  &\mi1 &1 & & 1&\mi1 & \dots&1 & \mi1&1 &\mi1 &1 &\mi1 & & \\
              n\ma3&  & & & \mi1& 1& &\dots &1 & \mi1&1 &\mi1 &1 &\mi1 & & \\
              \vdots&  & & &&& & \ddots & \ddots &&& \vdots  & \vdots & \vdots & \vdots & \vdots \\
              2n\mi3&  &  &  &  &  & && & \mi 1 & 1 &  & 1 & \mi1 &  & \\
              2n\mi2&  &  &  &  &  & && &&& \mi1 & 1 &  &  &
            \end{array}\right];\]
Now, since the determinant of $\B^1$ is 1 and $\C^2$ has been turned into a zero matrix, we can write $|\det(\A^1)|=|\det(\E^2)|$. To simplify $\E^2$, first add row $n\mi2$ to row $i$, where $1\leq i\leq n\mi3$. Second, add to row $n\ma i\ma 1$ the updated rows $j$, where $1\leq i \leq n\mi4$ and $1\leq j \leq n\mi i\mi 4$ and subtract from it row $n\mi2$. Third, add rows 1 to $n\mi4$ and rows $n\ma2$ to $2n\mi2$ to row $n\ma1$. Fourth, eliminate the nonzero entries of row $n$ from column $(1,2)$ to $(n\mi2, 1)$ by following similar steps that was used to eliminate the entries of row $(2n\mi2)$ in the proof of \cref{theorem:T4 FDI proof}. The resulting matrix $\E^3$ from the above elimination steps possesses the following form:

\[\lH \E^3=\scriptsize\arraycolsep=.25pt\def\arraystretch{1}\left[\begin{array}{lccccccccccccccc}
              \mcc{1}{i}& (1,\xlH2) & (2,\xlH1) & (1,\xlH3) & (3,\xlH1) & (1,\xlH4) & (4,\xlH1) & \dots & (1,\xlH n\mi3) & (n\mi3,\xlH1)& (1,\xlH n\mi2) & (n\mi2,\xlH1)&(1,\xlH n\mi1) & (n\mi1,\xlH1) & (1,\xlH n) & (n,1) \\
              1&  &  &  &  & & & & & & \mi1 & 1 & & & &  \\
              2&  &  &  &  & & & &\mi1 & 1 &  &  & & &  &  \\
              \vdots &  &  &  & &&& \iddots &\iddots&&&& \vdots  & \vdots&&\\
              n\mi4&  &  &\mi1 &1 &  &  &  &  & & & &  &  &  &  \\
              n\mi3& \mi1 & 1 & & &  &  &  &  & & & &  &  &  &  \\
              n\mi2&  & & & & & & & & & & & 1 & \mi1 &  &  \\
              n\mi1&  &  &  &  &  &  &  &  && &   && & & \mi1\\
              n&  &  &  &  &  &  & &  &  & & & \beta & \alpha & & \alpha \mi 1\\
              n\ma1&  & & & & & & & & & & & 1 & &\mi1 & \\
              n\ma2&  &\mi1 &1 & & &  & & & & & & & & & \\
              n\ma3&  & & & \mi1& 1& & & & & & & & & & \\
              \vdots&  & & &&& & \ddots & \ddots &&& \vdots  & \vdots & \vdots & \vdots & \vdots \\
              2n\mi3&  &  &  &  &  & && & \mi 1 & 1 &  & & &  & \\
              2n\mi2&  &  &  &  &  & && &&& \mi1 & 1 &  &  &
            \end{array}\right];\]
            
where, \[\alpha=\mi n\ma4;\mmmH \beta=2-\sum_{k=1}^{n-4}(k\mi1)=\frac{(n\mi1)\mi(n\mi3)(n\mi5)}{2};\]
Using the same reasoning as in the proof of \cref{theorem:T4 FDI proof}, it can be concluded that the non-singularity of $\E^3$ depends on the non-singularity of the following $4\xlH\times\xlH4$ matrix:
\[\left[\begin{array}{cccc}
    1     & \mi1         & \,\,\,0 & \,\,\,0 \\
    0     & \,\,\,0      & \,\,\,0 & \mi 1 \\
    \beta & \,\,\,\alpha & \,\,\,0 & \alpha\mi1\\
    1     & \,\,\,0      & \mi1    & \,\,\,0
  \end{array}\right].\]
The symbolic determinant of this matrix is $\alpha \ma \beta$, which equals 0 when,
\[\alpha\ma\beta=0\Leftrightarrow n=\frac{7}{2}\pm \frac{\sqrt{17}}{2}\notin \mathbb{Z}.\] 
which implies that the characteristic vectors produced by CP\Tsix\, are affinely independent, thereby establishing that \Tsix\,VI is facet defining for $n\geq4$. $_{\square}$\mmV
\end{pf} 
 
\subsection{\Tnine\,FDI Differences Matrix and Proof}\label{app:T9_proof}
To introduce the differences matrix for the \Tnine\,FDI proof the same relabeling of alternatives was used as in the previous subsection. After iteratively subtracting several rows of $\X$ as described in the below proof, the generating set for $\bar \X$ is given as: 
\[\arraycolsep=.75pt\def\arraystretch{0.9}\scriptsize
\begin{array}{|l|l|c||c|c|}\hline
j& \multicolumn{1}{|c|}{I_j:=i\in} &  p,q \T{-ranges}  & (k,\ell)\in \arcN: \bar{X}_{i,(k,\ell)} =1 & (k,\ell)\in \arcN: \bar{X}_{i,(k,\ell)} =\mi1 \\\hline
1 & \{1\} & - & (2, 1) & (n, 2) \\
2 & \{ 2n(p\mi1)\mi p^2\ma3 \}  & \ivl{1}{p}{n\mi2}\phantom{,p} & (n,\{1,p\ma1\})  & - \\
3 & \{ 2n(p\mi1)\mi p^2\ma4 \}  & \ivl{1}{p}{n\mi2}\phantom{,p} & - & (\{1,p\ma 1\}, n) \\
4 & \{ (2n\mi p)(p\mi1)\mi p\ma 2q\ma3\}& \ivl{1}{p}{n\mi3},\phantom{p}& (p\ma q\ma1,\{1,p\ma 1\}) & - \\
&&\ivl{1}{q}{n\mi2\mi p}&&\\
5 &\{ (2n\mi p)(p\mi1)\mi p\ma 2q\ma4\}& \ivl{1}{p}{n\mi3},\phantom{p} & - & (\{1,p\ma 1\}, p\ma q\ma 1)\\
&&\ivl{1}{q}{n\mi2\mi p}&&\\
6 & \{ 2np\mi (p\ma 1)^2+2 \} & \ivl{1}{p}{n\mi3}\phantom{,p} & (1,\{\ell \xlH>\xlH p\ma1\}); & (\{p\ma1\}\cup \{k \xlH>\xlH p\ma3\}, 1 );\\
& & & (p\ma2, n) & (n,p\ma2)\\
7 & \{ n(n\mi 2)\ma p \}  & \ivl{1}{p}{n\mi3}\phantom{,p} & (n\mi p\mi 1, \{1,n\mi 1\}) & (\{1,n\mi 1\}, n\mi p\mi 1) \\
8 & \{ n(n\mi 1) \mi 2 \}  & - & (n\mi1, n\mi 2) & (n\mi 2, n\mi 1) \\
9 & \{ n(n\mi 1) \mi 1 \}  & - & -  & (2,n) \\
10 & \{n(n\mi1)\} & \ivl{3}{p}{n\mi2}\phantom{,p} & (1, \{\ell \xlH> 1\xlH\}); (n,\{1<\ell \xlH<\xlH n\}); & - \\
&&& (2, \{\ell \xlH> 2 \xlH\});(p,\{p\xlH< \xlH \ell \xlH<\xlH n\}) & \\
\hline 
\end{array}\]

\begin{theorem}[\Tnine\,FDI]
\Tnine\,VI is an FDI of ${\bf P}^n_{WO}$, for any $n\ge4$.
\end{theorem} 
\begin{pf}
 As in the previous theorem, fix $i_1=1$ and $i_2 = n$ or assume that a corresponding relabeling of the alternatives is performed a priori. It is straightforward to verify that all points yielded by the modified version of CP\Tsix\,, denoted here as CP\Tnine\ for ease of explanation (the details of which can be found in page 21 of the main paper), belong to the ten ranking structures of \cref{T9_List} satisfy inequality \eqref{VI_T9} at equality. The rest of this proof follows almost the same steps as in the proof of \cref{theorem:T6 FDI proof}. Therefore we sketch below only minor differences and where a change in the structure occurs due to the differences between  CP\Tsix\, and CP\Tnine. First, most entries of the two difference matrices, both denoted by $\bar \X$, are the same except those in rows $1, n(n\mi1)\mi1$ and $n(n\mi1)$. The first row has only a $1$ under column $(2, 1)$ and a $\mi1$ under column $(n, 2)$, row $n(n\mi1)\mi1$ has a $\mi1$ under column $(2,n)$ with no other nonzero entry and finally row $n(n\mi1)$ has a binary structure that corresponds to a alternative-ordering where, item $i_1=1$ is in the first position, items $\{2,n\}$ are tied in the second position and the rest of the alternatives follow a lexicographical linear ordering. Second, due to the difference between the entries in the first row of $\A^0$ generated from CP\Tsix\, and CP\Tnine, the following changes are needed to convert $\A^0$ into $\A^1$: eliminate the entries in the first $(n\mi1)(n\mi2)$ columns of row $2n(i \mi1)\mi i^2\ma2$ by first adding the second row to row $1$ and then by adding rows $2n(i \mi1)\mi i^2\ma3$ and $2n(i \mi1)\mi i^2\ma4$ to row $2n(i \mi1)\mi i^2\ma2$, for $i=2,\dots,n\mi2$. Third, submatrix $\C^1$ has an additional $\mi1$ in row $n\mi1$ under column $(2,n)$, which can be eliminated by subtracting from it the second row of $\B^1$. The final change is related to difference in the structure of $\E^2$. The additional operations performed to eliminate the nonzero entry of $\C^1$ only affects two rows of $\E^2$; row $n\mi1$ which now consists of a $1$ under column $(1,n)$ and row $n\ma1$ which in addition to the $1$ under column $(2,1)$ has another $1$ under column $(n,1)$ instead of a $\mi1$ under column $(1,n)$. After the same elimination steps are performed on $\E^2$ as in the proof of \cref{theorem:T6 FDI proof} we get the following structure for $\E^3$:
\[\lH \E^3=\scriptsize\arraycolsep=.25pt\def\arraystretch{1}\left[\begin{array}{lccccccccccccccc}
              \mcc{1}{i}& (1,\xlH2) & (2,\xlH1) & (1,\xlH3) & (3,\xlH1) & (1,\xlH4) & (4,\xlH1) & \dots & (1,\xlH n\mi3) & (n\mi3,\xlH1)& (1,\xlH n\mi2) & (n\mi2,\xlH1)&(1,\xlH n\mi1) & (n\mi1,\xlH1) & (1,\xlH n) & (n,1) \\
              1&  &  &  &  & & & & & & \mi1 & 1 & & & &  \\
              2&  &  &  &  & & & &\mi1 & 1 &  &  & & &  &  \\
              \vdots &  &  &  & &&& \iddots &\iddots&&&& \vdots  & \vdots&&\\
              n\mi4&  &  &\mi1 &1 &  &  &  &  & & & &  &  &  &  \\
              n\mi3& \mi1 & 1 & & &  &  &  &  & & & &  &  &  &  \\
              n\mi2&  & & & & & & & & & & & 1 & \mi1 &  &  \\
              n\mi1&  &  &  &  &  &  &  &  && &   && & 1 &\\
              n&  &  &  &  &  &  & &  &  & & & \beta & \alpha & & \alpha \mi 1\\
              n\ma1&  & & & & & & & & & & & 1 & & & 1\\
              n\ma2&  &\mi1 &1 & & &  & & & & & & & & & \\
              n\ma3&  & & & \mi1& 1& & & & & & & & & & \\
              \vdots&  & & &&& & \ddots & \ddots &&& \vdots  & \vdots & \vdots & \vdots & \vdots \\
              2n\mi3&  &  &  &  &  & && & \mi 1 & 1 &  & & &  & \\
              2n\mi2&  &  &  &  &  & && &&& \mi1 & 1 &  &  &
            \end{array}\right];\]
where, $\alpha$ and $\beta$ have the same value as in the proof of \cref{theorem:T6 FDI proof}. Now, applying the same reasoning for the non-singularity of $E^3$ as in the proof of the previous theorem, it can be concluded that the determinant of this sub-matrix depends only on the following $4\xlH\times\xlH4$ matrix:
\[\left[\begin{array}{cccc}
    1     & \mi1         & \,\,\,0 & 0 \\
    0     & \,\,\,0      & \,\,\,1 & 0 \\
    \beta & \,\,\,\alpha & \,\,\,0 & \alpha\mi1\\
    1     & \,\,\,0      & \,\,\,0 & 1
  \end{array}\right].\]
The determinant of this matrix is $\beta\ma1$, which equals 0 when,
\[\beta\ma1=0\Leftrightarrow \frac{(n\mi 1)\mi(n\mi3)(n\mi5)}{2}\ma 1 = n^2\mi9n\ma14 = (n\mi2)(n\mi7)=0 \Leftrightarrow n=7\text{ or }2.\]
Therefore, for all $n\geq 4$, except $n=7$, the $n(n\mi1)$ characteristic vectors generated by CPT2-3 that satisfy inequality \eqref{VI_T9} at equality are linearly independent and, thus,  affinely independent. To prove that the result  holds for $n=7$ as well, subtract row $n(n\mi1)$ of matrix $\X$ from row $i$, where $i=1,2....n(n\mi1)\mi1$, to yield a new matrix $\hat \X$. It is straightforward to verify that the first $41$ vectors of $\hat \X$ are linearly independent (i.e., the row rank is 41) and, therefore, the $42$ characteristic vectors are indeed affinely independent. Hence, the result also holds for  $n=7$, and \Tnine\,VI is facet defining for any $n\geq4$. $_{\square}$\mmV
\end{pf} 

\bibliographystyle{cas-model2-names}

\bibliography{cas-refs}

\end{document}